\documentclass[sn-mathphys-num]{sn-jnl}% Math and Physical Sciences Numbered Reference Style
\usepackage{graphicx}%
\usepackage{multirow}%
\usepackage{amsmath,amssymb,amsfonts}%
\usepackage{amsthm}%
\usepackage{mathrsfs}%
\usepackage[title]{appendix}%
\usepackage{xcolor}%
\usepackage{textcomp}%
\usepackage{manyfoot}%
\usepackage{booktabs}%
\usepackage{algorithm}%
\usepackage{algorithmicx}%
\usepackage{algpseudocode}%
\usepackage{listings}%
\usepackage{cuted,multicol,lipsum}
\newtheorem{case}{Case}
\theoremstyle{thmstyleone}%
\newtheorem{theorem}{Theorem}%  meant for continuous numbers
\theoremstyle{thmstyletwo}%
\newtheorem{example}{Example}%
\newtheorem{remark}{Remark}%

\theoremstyle{thmstylethree}%

\begin{document}
\title[Article Title]{On the solutions of coupled nonlinear time-fractional diffusion-reaction system with time delays}

	%%%%%%--------------------------------------------------------------------
	%%%%%%          Title of the Paper and Acknowledgement
	%%%%%%--------------------------------------------------------------------
%	\title{On the solutions of {coupled nonlinear time-fractional diffusion-reaction system} with time delays
	%	\thanks{Acknowledgment: The first author (K.S.P.) is thankful for the financial support in the form of International Mathematical Union Breakout Graduate fellowship-2023 (IMU-BGF-2023-06) provided by the  IMU, Germany. Another author (M.L.) is supported by a Department of Science and Technology, India-SERB National Science Chair position (NSC/2020/000029).}
%	}
	%%%%%%--------------------------------------------------------------------
	%%%%%%         Authors,, Affiliations and email ids
	%%%%%%--------------------------------------------------------------------
	\author[1]{\fnm{K.S.} \sur{Priyendhu}}\email{ ks\_priyendhu@cb.students.amrita.edu}\equalcont{These authors contributed equally to this work.}
	
	\author*[1]{\fnm{P.} \sur{Prakash}}\email{vishnuindia89@gmail.com \,\&\, p\_prakash@cb.amrita.edu}
	\equalcont{These authors contributed equally to this work.}
	
	\author[2]{\fnm{M.} \sur{ Lakshmanan}}\email{lakshman.cnld@gmail.com}
	\equalcont{These authors contributed equally to this work.}
	
	\affil*[1]{\orgdiv{Department of Mathematics}, \orgname{Amrita School of Physical Sciences}, \orgaddress{\city{Coimbatore}, \street{Amrita Vishwa Vidyapeetham}, \country{INDIA}}}
	
	\affil[2]{\orgdiv{Department of Nonlinear Dynamics}, \orgname{Bharathidasan University},\orgaddress{  \city{Tiruchirappalli},\postcode{620 024}, \country{INDIA}}}
	
	%%%%%%--------------------------------------------------------------------
\abstract{
		In this article, we systematically explain how to apply the analytical technique called the invariant subspace method to find various types of analytical solutions for a  coupled  nonlinear time-fractional system of  partial differential equations  with time delays. Also, the present work explicitly studies a systematic way to obtain various kinds of finite-dimensional  {invariant  vector spaces} for the  coupled nonlinear time-fractional diffusion-reaction (DR) system  with time delays under the two {distinct} fractional derivatives, {namely} (a)  the Riemann-Liouville fractional partial time derivative  and (b)  the Caputo fractional  partial time  derivative.  Additionally, we provide  details of deriving analytical solutions in the generalized separable form for the initial and boundary value problems {(IBVPs)} of the     coupled nonlinear time-fractional DR system  with multiple time delays through the obtained   {invariant vector spaces} under the considered two time-fractional derivatives.
}

	
%	\begin{classification}
%		35R11; 35Bxx; 35-XX; 35Cxx
%	\end{classification}
	\maketitle
	
	\section{Introduction}\label{sec1}
In recent years, applications of fractional differential equations in  different fields of science and engineering {have expanded} tremendously due to the {attractive features} of fractional-order derivatives  in describing memory and hereditary characteristics when compared to integer-order {derivatives} \cite{Diethelm2010,Tarasov2011}.
   One of the prominent classes of arbitrary-order (fractional-order) derivatives is those with weakly singular kernels, as they possess peculiar properties like the nonlocal nature, irreversibility, fractal properties, and long-term memory behavior \cite{Diethelm2010,Tarasov2011}. % These properties contribute to the accuracy and precision of complex and unpredictable processes. The anomalous diffusion process is a classic example of an unpredictable physical process, which can be successfully modelled using arbitrary-order derivatives  The most straightforward mathematical formulation of the anomalous diffusion process \cite{bakkyaraj2015}, with the help of the fractional-order derivatives and integrals,  can be written as
%	\begin{equation}\label{diffusion}
%		\partial^\alpha_tu(x,t)=a \partial_{xx} u(x,t),0<\alpha\leq2,
%	\end{equation}
%where $\partial^{\alpha}_{t}u(x,t)$ is the $\alpha$-th order {fractional partial time derivative} of  $u(x,t),$ $\partial_{xx}u(x,t)=\dfrac{\partial^2 u(x,t)}{\partial x^2},$ $x\in\mathbb{R},t>0,$ and $a$ is a constant that depends on the diffusion nature of the substance considered.
%Note that the range of  $\alpha$ determines the nature of the physical process that the equation \eqref{diffusion} represents.
%Note the range of $\alpha$ determine the nature of the physical process that the equations \eqref{diffusion} represents.
%If $0<\alpha<1,$  equation \eqref{diffusion} corresponds to sub-diffusion phenomena \cite{bakkyaraj2015}.  Also, equation   \eqref{diffusion} represents  normal diffusion and the wave phenomena when $\alpha=1 $ and $\alpha=2,$ respectively. Moreover, equation \eqref{diffusion} is called a diffusion-wave   phenomenon \cite{bakkyaraj2015} if $1<\alpha<2$.

We know that the solutions of fractional PDEs (FPDEs) are essential because they will help us to understand the physical processes and phenomena of the considered systems.  However, studying analytical solutions for nonlinear FPDEs is difficult since the fractional-order derivatives with weakly singular kernels fail to follow crucial standard rules such as the Leibniz rule, the semigroup property, and the chain rule \cite{ Diethelm2010,Tarasov2011}.
 %Due to violating standard rules of fractional-order derivatives, all classical methods of integer-order PDEs cannot be extended to solve FPDEs.
  Recent studies show that  the Lie symmetry analysis \cite{bakkyaraj2015,prakash2021}, %the Adomian decomposition method \cite{gejji2005,jafari2006},
 {the  method of variable separation} \cite{ru2020,rui2022,ruinonlineardyn2022}, and the invariant subspace method \cite{gazizov2013,harris2013,sahadevan2015,choudary2017,harris2017,rui2018,sahadevan2017a,prakash2020a,prakash2020b,gara2021,prakash2022a,prakash2022b,priyendhu2023,prakash2023,prakash2024} are effective {mathematical} methods to find analytical solutions for various kinds of nonlinear FPDEs.

	The  invariant subspace method for {obtaining} the analytical solution of nonlinear integer-order (classical) PDEs was  initially discussed by Galaktionov and Svirshchevskii \cite{Galaktionov2007}. This method was further developed and   investigated by  Ma \cite{ma2012a} and others \cite{ma2012b,ma2012c,qu2009, Polyanin2024,polyanin2014coupled}.
 %In \cite{gazizov2013}, Gazizov and Kasatkin were the first {to solve FPDEs through the  invariant subspace method.}
 Recently, {a number of} studies have been carried out to develop this method for several classes of nonlinear FPDEs in the literature \cite{gazizov2013,harris2013,sahadevan2015,choudary2017,harris2017,rui2018,sahadevan2017a,prakash2020a,prakash2020b,gara2021,prakash2022a,prakash2022b,prakash2023,priyendhu2023,prakash2024}.
However,  only a very few {of these} studies have investigated the efficacy of the invariant subspace method in finding the generalized separable analytical solutions for nonlinear FPDEs with time delays \cite{prakash2020b,prakash2022a,priyendhu2023}.
  {Moreover, both local and nonlocal integrable equations with classical derivatives have been studied extensively using different methods such as Darboux	transformations \cite{paper1}   or
 	Riemann-Hilbert transforms \cite{paper2}. Also, we would like to mention that multi-component coupled integrable  equations  have been solved through the symmetry point of view as discussed in \cite{paper3,paper4}. In the future, we believe that the theory of the invariant subspace method can also be extended to solve such classes of local and nonlocal multi-component coupled integrable equations with both integer and fractional-order derivatives. }

In the literature, systems of PDEs with multiple time delays  have been studied to mathematically describe different physical processes that arise in the field of  science and engineering, such as neurology, economics, image processing and so on \cite{Lakshmanan2010,Polyanin2024}.
 In this work, we aim to  systematically apply the invariant subspace method to find various types of analytical solutions for {coupled nonlinear time-fractional systems of PDEs with time delays} for the first time.
 %We give the systematic details of  constructing invariant vector spaces for the considered {coupled nonlinear time-fractional system} with time delays involving nonlinear or linear delay terms.
 Additionally, we present the effectiveness of  the invariant subspace method to {obtain} solutions for a {coupled nonlinear time-fractional DR system} with time delays in the following form
\begin{eqnarray}
	\begin{aligned}\label{DREs}
		&	\left( 	{\partial^{\alpha_1}_{t}} u_1,
		{\partial^{\alpha_2}_{t}} u_2
		\right)	= \mathbf{R}[\mathbf{U},\hat{\mathbf{U}}],0<\alpha_1,\alpha_2\leq2,
	\end{aligned}
\end{eqnarray}
 $\mathbf{R}[\mathbf{U},\hat{\mathbf{U}}]
=
\left(
\mathcal{R}_1(\mathbf{U},\hat{\mathbf{U}}),
\mathcal{R}_2(\mathbf{U},\hat{\mathbf{U}})
\right)$ with
\begin{eqnarray}
	\begin{aligned}\label{components of dre}
		\mathcal{R}_k(\mathbf{U},\hat{\mathbf{U}})=&\partial_x[(a_{k2}u_2+a_{k1}u_1+a_{k0})\partial_x u_1
		+(b_{k2}u_2+b_{k1}u_1+b_{k0})\partial_x u_2 ]
		\\&	
		+ \gamma_{k5}u_1u_2	+\gamma_{k4}u_2^2+\gamma_{k3}u_1^2 +\gamma_{k2}u_2 	 +\gamma_{k1}u_1
		+\gamma_{k0} +
		\hat{\gamma}_{k1}\hat{u}_1 +\hat{\gamma}_{k2}\hat{u}_2,
	 	\end{aligned}
\end{eqnarray}
{where}  ${\partial^{\alpha_k}_{t}}(\cdot)$ denotes the $\alpha_k$-th order fractional partial time derivative of $u_k(x,t)$  in the sense of   {either}   Caputo  ${^C\partial^{\alpha_k}_{t}}(\cdot),$ {or}  Riemann-Liouville   ${^{RL}\partial^{\alpha_k}_{t}}(\cdot),$  $\hat{u}_k=\hat{u}_k(x,t-\tau_k),$  { $\tau_k>0$ represents the  time delay of the $k$-th component,}   $x\in\mathbb{R} $ is the space variable, $t>0 $ is the time variable,  and $a_{kj},b_{kj},\gamma_{kl},\hat{\gamma}_{k1},\hat{\gamma}_{k2}\in\mathbb{R},j=0,1,2,l=0,1,\dots,5,k=1,2.$
  We know that the Caputo {fractional  partial time derivative}  ${^C\partial^{\alpha}_{t}}(\cdot)$ and  the Riemann-Liouville {fractional  partial time derivative} ${^{RL}\partial^{\alpha}_{t}}(\cdot)$ of order $\alpha>0$ are defined in \cite{Diethelm2010} as follows:
 \small{ \begin{eqnarray}
 		\label{Caputo}	{^C\partial^{\alpha}_{t}}f(x,t)= \mathit{I}_t^{n-\alpha} \left[ \partial^n_tf(x,t)\right],n-1<\alpha\leq n,
 		\\
 		\label{RL}		{^{RL}\partial^{\alpha}_{t}}f(x,t)= \partial^n_t\left[\mathit{I}_t^{n-\alpha}  f(x,t)\right] , n-1<\alpha\leq n,
 \end{eqnarray}}
 respectively. Note that  $\mathit{I}_t^{n-\alpha}f(x,t)=
 \dfrac{1}{\Gamma(n-\alpha)}\int\limits_{0}^t(t-v)^{n-\alpha-1}f(x,v)dv$ is the Riemann-Liouville
 time-fractional integral of order $n-\alpha>0,  $  where $\alpha\in(n-1,n],$ $n\in\mathbb{N}.$
In \cite{Polyanin2024},
Polyanin and Zhurov have shown the applicability of the invariant subspace method to {obtain} the generalized separable analytical solutions for some special cases of  \eqref{DREs}   when $\alpha_k=1,$ and $a_{ki}=b_{ki}=0,i,k=1,2$. For $\alpha_k=1,$    the above system \eqref{DREs}  has  potential applications in population dynamics models like the Lotka-Volterra system, the Belousov–Zhabotinsky oscillating reaction models,   and many others \cite{Polyanin2024,polyanin2014coupled}.

The main objective of this work is to study the analytical solutions for the {IBVPs} of the  given  {coupled nonlinear time-fractional DR system} \eqref{DREs}  with time delays  under the two familiar fractional {partial time derivatives} that are the Caputo  {fractional partial time derivative} ${^C\partial^{\alpha_k}_{t}}(\cdot),$ and the Riemann-Liouville {fractional  partial time derivative} ${^{RL}\partial^{\alpha_k}_{t}}(\cdot)$ of order $\alpha_k,0<\alpha_k\leq2,k=1,2.$
Additionally, we wish to point out that this work studies the given {coupled nonlinear time-fractional DR system} \eqref{DREs} under two different {fractional  partial time derivatives} through the invariant subspace method for the first time in the literature.
      In this {paper}, we consider the given {coupled nonlinear time-fractional DR system}  \eqref{DREs} with the
   initial  and the  boundary conditions as given below.
   \begin{itemize}
   	\item[$\bullet$] For the Caputo partial time derivative ${\partial^{\alpha}_{t}}(\cdot)={^C\partial^{\alpha}_{t}}(\cdot)$, we consider
   	\begin{eqnarray}
   		\begin{array}{ll}
   			&\label{ic-bc caputo}u_k(x,0)=  {G}_{k,0}(x),
   			\partial_{t}[u_k(x,t)]\Big{|}_{t=0}= {G}_{k,1}(x) ,
   				\\&		u_k(0,t)=F_{k,0}(t),
   	   			  			\text{ and }
   		u_k(\lambda,t)=F_{k,\lambda}(t),k=1,2,\lambda\in\mathbb{R}.
   		\end{array}
   	\end{eqnarray}
      	\item[$\bullet$] For the Riemann-Liouville partial time derivative ${\partial^{\alpha}_{t}}(\cdot)={^{RL}\partial^{\alpha}_{t}}(\cdot),$  we consider
   	\begin{eqnarray}
   		\begin{array}{ll}
   			&\label{ic-bc rl}
   			{^{RL}\partial}^{\alpha_k-l}_{t}[u_k(x,t)]\Big{|}_{t=0}= \xi_{k,l}(x) ,
   				u_k(0,t)=f_{k,0}(t),
   		\\&	\text{ and }
   			u_k(\lambda,t)=f_{k,\lambda}(t),k,l=1,2,\lambda\in\mathbb{R}.
   			\end{array}
   	\end{eqnarray}
   \end{itemize}

 This article is organized as follows: Section \ref{sec2} presents the systematic
 development of the theory of invariant subspace method to derive solutions in generalized separable form for a  {coupled nonlinear time-fractional system of   PDEs} with time delays.
 %Also, we provide the details of how the invariant vector spaces of the {coupled nonlinear time-fractional systems of PDEs} with time delays help to find their solutions in generalized separable form.
 In Section \ref{sec3}, we present the applicability of the invariant subspace method to derive solutions for the IBVPs of the given system \eqref{DREs}  under the Riemann-Liouville and Caputo {fractional partial time derivatives. } {Also, we provide some extensions of the invariant subspace method to obtain analytical solutions for the {$2$-component and multiple-component coupled nonlinear time-fractional  systems of  PDEs} with several linear time delays in Section \ref{several}. }
  Finally, we present the conclusion of the obtained results for the discussed method in Section  \ref{sec4}.
\section{Invariant vector spaces and  solutions of  coupled nonlinear time-fractional system of   PDEs with time delays}\label{sec2}
This section explains the algorithmic way to find generalized separable analytical solutions for the  {coupled nonlinear time-fractional system of   PDEs} with time delays  using the invariant subspace method.
Thus, we consider the {coupled nonlinear time-fractional system of   PDEs} with time delays involving nonlinear terms of  the form  \begin{eqnarray}
	\begin{aligned}\label{1+1:delay}
	&	\left( 	{\partial^{\alpha_1}_{t}} u_1,
		{\partial^{\alpha_2}_{t}} u_2
		\right)
		={\mathbf{{{M}}}}[\mathbf{U},\hat{\mathbf{U}}],
			%	\left(		{\mathcal{{{M}}}_1}[\mathbf{U},\hat{\mathbf{U}}],		{\mathcal{{{M}}}_2}[\mathbf{U},\hat{\mathbf{U}}]		\right) ,
		\alpha_1,\alpha_2>0,
		\\&	u_k(x,t)=\phi_k(x,t),t\in[-\tau_k,0],k=1,2,
	\end{aligned}
\end{eqnarray}
where $	{\partial^{\alpha_k}_{t}}(\cdot)$ represents an  $\alpha_k$-th order  {fractional partial time derivative}  of  $\mathbf{U}=(u_1,u_2)$ in the sense of either Caputo ${^C\partial^{\alpha_k}_{t}}(\cdot)$ given in \eqref{Caputo} or the Riemann-Liouville  ${^{RL}\partial^{\alpha_k}_{t}}(\cdot)$ given in \eqref{RL},
$\hat{\mathbf{U}}=(\hat{u}_1,\hat{u}_2),$
$ u_k=u_k(x,t),$
$\hat{u}_k=\hat{u}_k(x,t-\tau_k),$
$x\in\mathbb{R}, \tau_k,t>0,$
\begin{eqnarray*}
	\begin{aligned}
{\mathbf{{{M}}}}[\mathbf{U},\hat{\mathbf{U}}]=&	\left(		{\mathcal{{{M}}}_1}[\mathbf{U},\hat{\mathbf{U}}],		{\mathcal{{{M}}}_2}[\mathbf{U},\hat{\mathbf{U}}]		\right)
\text{  and  }
\\	{\mathcal{{{M}}}_k}[\mathbf{U},\hat{\mathbf{U}}]
	=&	{\mathcal{{{M}}}_k}\big( x,
	u_1,\hat{u}_1,
		\partial_x u_1,
			\partial_x \hat{u}_1, 	\dots,
			\partial_x^{m_k} u_1,	\partial_x^{m_k} \hat{u}_1,
u_2,
\hat{u}_2,
\partial_x u_2,
%\partial_x \hat{u}_2,
\dots,	\partial_x^{m_k}  u_2,
\partial_x^{m_k} \hat{u}_2
	\big)
\end{aligned}
\end{eqnarray*}
is the function corresponding to the nonlinear delay differential operator of order $m_k,m_k\in\mathbb{N},$ with $ \partial_x^{j}(\cdot)=\dfrac{\partial^j}{\partial x^j}(\cdot),j=1,2,\dots,m_k,$ $k=1,2.$ Hereafter, we will call the function ${\mathcal{{{M}}}_k}[\mathbf{U},\hat{\mathbf{U}}]$ as the operator for clarity.

Now, we recall the definition of the Cartesian product of two vector spaces.  So, let  us consider the finite-dimensional vector spaces $Y_{k,n_k}$  with  dim$(Y_{k,n_k})=n_k,k=1,2.$
Then, the Cartesian product of two vector spaces $Y_{k,n_k},k=1,2$  is defined as
\begin{eqnarray}
	\begin{aligned}
	\mathbf{Y}_n:=&Y_{1,n_1}\times
	Y_{2,n_2}
	=	\{ y=(y_1,y_2)|
	y_k\in Y_{k,n_k},k=1,2\}
	\end{aligned}
	\end{eqnarray}
along with the  vector addition $y+z:=(y_1,y_2)+(z_1,z_2)=(y_1+z_1,y_2+z_2),$
and scalar multiplication
$cy:=c(y_1,y_2)=(cy_1,cy_2),$ for all $y=(y_1,y_2),z=(z_1,z_2)\in 	\mathbf{Y}_n,$ and $ c\in\mathbb{R}.$

The following theorem explains the dimension of the  Cartesian product of two vector spaces $\mathbf{Y}_n=Y_{1,n_1}\times
Y_{2,n_2}.$
\begin{theorem}
	\label{axlertheorem}\cite{Axler2014}
	If the dimension of the vector spaces $Y_{k,n_k}$ is  $n_k$ $ (=\text{dim}(Y_{k,n_k})),$ then the dimension of the Cartesian product $\mathbf{Y}_n$ of two vector spaces $Y_{k,n_k},k=1,2$  is
$n=n_1+n_2\, (:=	\text{dim}(\mathbf{Y}_n)).$
\end{theorem}
Now, the following steps provide details on how to develop the idea of the invariant subspace method for the given {coupled nonlinear time-fractional system of PDEs} \eqref{1+1:delay} with time delays.
\begin{enumerate}
	\item \label{step1}
	First, we should find the  invariant vector spaces for the given coupled nonlinear time-fractional system \eqref{1+1:delay} with time delays  using the solutions of linear homogeneous integer-order ODEs.
	\item \label{step2}
	Based on the obtained invariant vector spaces for the given coupled system \eqref{1+1:delay}, we can expect the solution as a generalized separable form for $(u_1,u_2)$ of the given {{system}} \eqref{1+1:delay} and reduce it into a system of fractional-order time derivative delay ODEs.
	\item \label{step3}
	Using   suitable analytical methods, we can solve the obtained system of fractional-order time derivative delay ODEs, from which we get the generalized separable analytical solution for  \eqref{1+1:delay}.
\end{enumerate}
Next, we systematically explain the steps mentioned above one by one for obtaining analytical solutions as a generalized separable type for the {present system} \eqref{1+1:delay}.
\subsection{Invariant vector spaces for ${\mathbf{{{M}}}}[\mathbf{U},\hat{\mathbf{U}}]$ given in \eqref{1+1:delay}}
Let us consider the system of linear homogeneous integer-order
ODEs {as follows:}
\begin{eqnarray}
	\begin{aligned}\label{lhodes}
		H_k[\lambda(x)]	:=&{\mathbf{D}^{n_k}_x\lambda(x)}+\mu_{k,n_k-1}{\mathbf{D}^{n_k-1}_x\lambda(x)}+\dots
		 +\mu_{k,0}\lambda(x)=0, k=1,2,
	\end{aligned}
\end{eqnarray}
where $\mathbf{D}^{i}_x=\dfrac{d^i}{d x^i},j=1,2,\dots,n_k, $ $n_k\in\mathbb{N},  $ represent the ordinary integer-order derivative operator, and $ \mu_{k,i},i=1,2,\dots,n_k-1,k=1,2$ are arbitrary constants. Now, we assume that the set of  $n_k$-linearly independent functions $\{\psi_{k,i}(x)|i=1,\dots,n_k\}$ form a basis (fundamental) set of  solutions for the linear homogeneous
integer-order ODEs $H_k[\lambda(x)]=0,k=1,2.$ From this, we can define  the component vector spaces (or solution spaces) $Y_{k,n_k}$ as the linear span of  $n_k$-linearly independent functions $\psi_{k,i}(x),i=1,2,\dots,n_k,$ for $k=1,2.$ That is,
\begin{eqnarray}\label{component spaces}
	\begin{aligned}
		Y_{k,n_k}=&\mathcal{L}
		\{\psi_{k,1}(x),\psi_{k,2}(x),\dots,\psi_{k,n_k}(x)\},
	%	\\=&\{ \omega_k(x)=\sum\limits_{i=1}^{n_k}\delta_{k,i} \psi_{k,i}(x)\ | \ \delta_{k,i}\in\mathbb{R},\\& H_k[\omega_k(x)]=0,  i=1,2,\dots,n_k\},
	\end{aligned}
\end{eqnarray}
where $\mathcal{L} $ denotes the linear span. % and $H_k[\cdot],k=1,2,$ are  given in \eqref{lhodes}.
Hence, the Cartesian product space $\mathbf{Y}_n=Y_{1,n_1}\times
Y_{2,n_2}$ can be read as follows:
\begin{eqnarray}
	\begin{aligned}\label{vector space}
		\mathbf{Y}_n=& \mathcal{L}\{ (\psi_{1,1}(x),0),\dots,  (\psi_{1,n_1}(x),0),
		  (0,\psi_{2,1}(x)),\dots,(0,\psi_{2,n_2}(x))\}
		\\=&
		\{\omega(x)=(\omega_1(x),\omega_2(x)) |\, \omega_k(x) \in Y_{k,n_k},  \omega_k(x)=\sum\limits_{i=1}^{n_k}\delta_{k,i} \psi_{k,i}(x), \delta_{k,i}\in\mathbb{R}, 	\\&	H_k[\omega_k(x)]=0,
		 i=1,2,\dots,n_k,k=1,2
		\}.
	\end{aligned}
\end{eqnarray}
We call the Cartesian product space $\mathbf{Y}_n$ is an invariant vector space for the vector nonlinear delay differential operator ${\mathbf{{{M}}}}[\mathbf{U},\hat{\mathbf{U}}]$ given in \eqref{1+1:delay} if any one of the following  invariant conditions holds:
\begin{eqnarray}
	\begin{aligned}
		&{\mathbf{{{M}}}}[\mathbf{U},\hat{\mathbf{U}}]\in \mathbf{Y}_n \, \forall \, \mathbf{U} \in \mathbf{Y}_n,
%	\\&	\text{ which means that }
	%	{\mathbf{{{M}}}}[\mathbf{Y}_n,\mathbf{Y}_n] \subseteq\mathbf{Y}_n,
		  \text{ or } &
	%	{\mathcal{{{M}}}_k}[\mathbf{U},\hat{\mathbf{U}}]\in  Y_{k,n_k}, \forall \, \mathbf{U} \in \mathbf{Y}_n,
%	\\&	\text{ that is }
		{\mathcal{{{M}}}_k}[\mathbf{Y}_n,\mathbf{Y}_n] \subseteq Y_{k,n_k},k=1,2.
	\end{aligned}
\end{eqnarray}
Thus, we can say that whenever  $u_k$ and $\hat{u}_k$ are obtained as $u_k=\sum\limits_{i=1}^{n_k}\delta_{k,i}\psi_{k,i}(x)$ and $
\hat{u}_k=\sum\limits_{i=1}^{n_k}\hat{\delta}_{k,i}
\psi_{k,i}(x),$ then  there exist functions  $\Lambda_{k,i},i=1,2,\dots,n_k,k=1,2,$  such that
	\begin{eqnarray}\label{invariant condition}
	\begin{aligned}
		{\mathcal{{{M}}}_k} [  ( u_1,u_2 ), ( \hat{u}_1,\hat{u}_2 )  ]
		=&
		\sum\limits_{i=1}^{n_k}
		\Lambda_{k,i}\Big(
		\delta_{1,1},
		\hat{\delta}_{1,1},
		\delta_{2,1},
			\ldots,
		\hat{\delta}_{1,n_1},
			\delta_{2,n_2},
		\hat{\delta}_{2,n_2}\Big)\psi_{k,i}(x).
	\end{aligned}
\end{eqnarray}
Next,   we  discuss how  to obtain different kinds of analytical solutions   using the invariant vector space $\mathbf{Y}_n$ for the  given {coupled nonlinear time-fractional system} \eqref{1+1:delay}  in the following generalized separable form
\begin{eqnarray}
	\mathbf{U}(x,t)=\big( \sum\limits_{i=1}^{n_1}\delta_{1,i}(t)\psi_{1,i}(x),\sum\limits_{i=1}^{n_2}\delta_{2,i}(t)\psi_{2,i}(x)\big) ,\label{gses form}
\end{eqnarray}
where $ \delta_{k,i}(t)$ are unknown functions of time variable $t$ that can be determined from the  given {coupled nonlinear time-fractional system} \eqref{1+1:delay} and the obtained space $\mathbf{Y}_n.$ The details are presented  in the following subsection.
\subsection{Solutions of the  given system \eqref{1+1:delay}}
In this subsection, we provide the necessary criteria to admit the required form of analytical solution \eqref{gses form} for the given {coupled nonlinear time-fractional system of PDEs}  \eqref{1+1:delay} with time delays, which is the primary key result of our study that is discussed in the following theorem.
%%%%%%%%%%%%%%%%%%%%%%%THEOREM%%%%%%%%%%%%%%%%%%%%%%%%%%%%
\begin{theorem}
	\label{main result}
	Let the vector nonlinear delay differential operator  ${\mathbf{{{M}}}}[\mathbf{U},\hat{\mathbf{U}}]$ given in \eqref{1+1:delay}  admit the invariant vector space
	%\begin{eqnarray}
	%	\begin{aligned}
	$		\mathbf{Y}_n= Y_{1,n_1}\times
			Y_{2,n_2}.
		%	\\=& \mathcal{L}\{ (\psi_{1,1}(x),0),\dots,  (\psi_{1,n_1}(x),0),
		%	\\&(0,\psi_{2,1}(x)),
		%	  \dots,(0,\psi_{2,n_2}(x))\}.
	$%	\end{aligned}
	%\end{eqnarray}
	Then the analytical solution of the given {coupled nonlinear time-fractional system}  \eqref{1+1:delay} is obtained as
	 \begin{eqnarray}	\begin{aligned}
		 	\mathbf{U}(x,t):=&(u_1(x,t),u_2(x,t))
 =\left( \sum\limits_{i=1}^{n_1}\delta_{1,i}(t)\psi_{1,i}(x),\sum\limits_{i=1}^{n_2}\delta_{2,i}(t)\psi_{2,i}(x)\right).\label{thm:gses form}
\end{aligned}	\end{eqnarray}
	Moreover, the functions $\delta_{k,i}(t),i=1,2,\dots,n_k,k=1,2$  are determined from the following system of fractional-order time derivative delay ODEs
	\begin{eqnarray}	\begin{aligned}\label{thm:sysodes}
			D^{\alpha_k}_t{\delta_{k,i}(t)}=&
			\Lambda_{k,i}\Big(
			\delta_{1,1}(t),
			\hat{\delta}_{1,1}(t-\tau_1),
			\delta_{2,1}(t),
		\hat{\delta}_{2,1}(t-\tau_2),
			\ldots, 	\hat{\delta}_{1,n_1}(t-\tau_1),
			\\&
					\delta_{1,n_1}(t),		\delta_{2,n_2}(t),
			\hat{\delta}_{2,n_2}(t-\tau_2)\Big), {i=1,2,\dots,n_k,k=1,2,}
		\end{aligned}
	\end{eqnarray}
where $ D^{\alpha_k}_t(\cdot),k=1,2,$ denotes the $\alpha_k$-th order fractional ordinary time derivative
%	is  an analytical solution for the given  {coupled nonlinear time-fractional system of PDEs} \eqref{1+1:delay} of PDEs with time delays 	$$\left( 	{\partial^{\alpha_1}_{t}} u_1,	{\partial^{\alpha_2}_{t}} u_2	\right) 	={\mathbf{{{M}}}}[\mathbf{U},\hat{\mathbf{U}}],\alpha_1,\alpha_2>0.$$
 in the sense of Caputo  {or} Riemann-Liouville.
\end{theorem}
\begin{proof}
	Suppose that $	\mathbf{Y}_n=Y_{1,n_1}\times
	Y_{2,n_2}$ is invariant under the vector nonlinear delay differential operator  $${\mathbf{{{M}}}}[\mathbf{U},\hat{\mathbf{U}}]=	\left(
	{\mathcal{{{M}}}_1}[\mathbf{U},\hat{\mathbf{U}}],
	{\mathcal{{{M}}}_2}[\mathbf{U},\hat{\mathbf{U}}]
	\right),$$ which is given in \eqref{1+1:delay}.
Thus, we get  \eqref{invariant condition}.
Now, due to the invariant vector space $	\mathbf{Y}_n $    of ${\mathbf{{{M}}}}[\mathbf{U},\hat{\mathbf{U}}],$  we assume that $u_k(x,t)=\sum\limits_{i=1}^{n_k}\delta_{k,i}(t)\psi_{k,i}(x)$ and $u_k(x,t-\tau_k)=\sum\limits_{i=1}^{n_k}\hat{\delta}_{k,i}(t-\tau_k)\psi_{k,i}(x)$. Hence, by using \eqref{invariant condition}, we obtain
%that ${\mathcal{{{M}}}_k}[\mathbf{U},\hat{\mathbf{U}}],k=1,2,$ satisfy the conditions as follows:
	\begin{eqnarray}\label{thm:invariant condition}
		\begin{aligned}
				{\mathcal{{{M}}}_k}[\mathbf{U},\hat{\mathbf{U}}]&=
			\sum\limits_{i=1}^{n_k}
			\Lambda_{k,i}\Big(
			\delta_{1,1}(t),
			\hat{\delta}_{1,1}(t-\tau_1),
			\delta_{2,1}(t),
	\hat{\delta}_{2,1}(t-\tau_2),
			\ldots,
			 			\delta_{1,n_1}(t)
			,
		\hat{\delta}_{1,n_1}(t-\tau_1),	
		\\&	\delta_{2,n_2}(t),
			\hat{\delta}_{2,n_2}(t-\tau_2)\Big)\psi_{k,i}(x),
		\end{aligned}
	\end{eqnarray}	{where $\Lambda_{k,i}(\cdot),i=1,2,\dots,n_k,k=1,2,$ are the expansion coefficients determined from the basis set of  $	\mathbf{Y}_n.$}
	Also, we know that  the Caputo ${^C\partial^{\alpha_k}_{t}}(\cdot)$ and the Riemann-Liouville ${^{RL}\partial^{\alpha_k}_{t}}(\cdot)$  {fractional partial time derivatives}  are linear differential operators. Then,  we get
	\begin{eqnarray}\label{thm:sysy}
		\begin{aligned}
			{\partial^{\alpha_k}_{t}}u_k:=&	{\partial^{\alpha_k}_{t}}\left( \sum\limits_{i=1}^{n_k}\delta_{k,i}(t)\psi_{k,i}(x)
			\right)
		= \sum\limits_{i=1}^{n_k}	D^{\alpha_k}_t\left[\delta_{k,i}(t)\right] \psi_{k,i}(x).
		\end{aligned}
	\end{eqnarray}
if  $u_k(x,t)=\sum\limits_{i=1}^{n_k}\delta_{k,i}(t)\psi_{k,i}(x),k=1,2.$

Since $\mathbf{Y}_n$ is an invariant vector space for  ${\mathbf{{{M}}}}[\mathbf{U},\hat{\mathbf{U}}] ,$  %$u_k(x,t)=\sum\limits_{i=1}^{n_k}\delta_{k,i}(t)\psi_{k,i}(x),$ and $u_k(x,t-\tau_k)=\sum\limits_{i=1}^{n_k}\hat{\delta}_{k,i}(t-\tau_k)\psi_{k,i}(x),k=1,2,$
 from \eqref{thm:invariant condition} and \eqref{thm:sysy},  we get
	$$
			\left( 	{\partial^{\alpha_1}_{t}} u_1,
			{\partial^{\alpha_2}_{t}} u_2
			\right)
			 ={\mathbf{{{M}}}}[\mathbf{U},\hat{\mathbf{U}}]
		$$
		{ admits the solution \eqref{thm:gses form} if }
	$$
	\begin{aligned}
		&	\sum\limits_{i=1}^{n_k}\Big[ D^{\alpha_k}_t[\delta_{k,i}(t)]-
		\Lambda_{k,i}\Big(
		\delta_{1,1}(t),
		\hat{\delta}_{1,1}(t-\tau_1),
		\delta_{2,1}(t),
		\hat{\delta}_{2,1}(t-\tau_2),
		\ldots,
		\delta_{1,n_1}(t)
		,
		\hat{\delta}_{1,n_1}(t-\tau_1),
	\\
	&	\delta_{2,n_2}(t),
			\hat{\delta}_{2,n_2}(t-\tau_2)\Big)\Big] \psi_{k,i}(x)=0.
	\end{aligned}
	$$
	Hence,  we get $u_k(x,t)=\sum\limits_{i=1}^{n_k}\delta_{k,i}(t)\psi_{k,i}(x),k=1,2,$  is  an analytical solution of the given {coupled nonlinear time-fractional system} \eqref{1+1:delay} with  the functions $\delta_{k,i}(t)$ {to be} determined by the given system of fractional-order time derivative delay ODEs \eqref{thm:sysodes}
 	whenever the given system \eqref{1+1:delay} admits $	\mathbf{Y}_n.$
\end{proof}
\subsection{Invariant vector spaces for a special case of ${\mathbf{{{M}}}}[\mathbf{U},\hat{\mathbf{U}}]$ given in \eqref{1+1:delay} }
This subsection presents the invariant vector spaces for the special  case of the given {coupled nonlinear time-fractional system of PDEs} \eqref{1+1:delay} with time delays, in which time delay occurs only in linear terms of $\hat{u}_1 $ and $\hat{u}_2. $

Let us consider the special case of the given {coupled nonlinear time-fractional system}  \eqref{1+1:delay} with linear time delay  terms  of $\hat{u}_1 $ and $\hat{u}_2 .$  Thus,  {{the vector nonlinear delay differential operator}} ${\mathbf{{{M}}}}[\mathbf{U},\hat{\mathbf{U}}]$  takes the following form,
\begin{eqnarray}
	\begin{aligned}\label{linear-delay operator}
		{\mathbf{{{M}}}}[\mathbf{U},\hat{\mathbf{U}}]
		=&
		\left(
		{\mathcal{{{M}}}_1}[\mathbf{U},\hat{\mathbf{U}}],
		{\mathcal{{{M}}}_2}[\mathbf{U},\hat{\mathbf{U}}]
		\right)
		\,	\text{along with }
		\\
		{\mathcal{{{M}}}_k}[\mathbf{U},\hat{\mathbf{U}}]
		=&	{\mathcal{{{M}}}_k} ( x,
		u_1,\partial_x u_1,	\dots,\partial_x^{m_k} u_1,
		u_2,
		\partial_x u_2,
		\dots,	\partial_x^{m_k}  u_2
		)
		+
		\eta_{k,1}	\hat{u}_1+
		\eta_{k,2}	\hat{u}_2,
	\end{aligned}
\end{eqnarray}
where $\eta_{k,i}\in\mathbb{R},k,i=1,2.$ Note that  $\mathbf{Y}_n$ is invariant under the above-vector nonlinear delay differential operator ${\mathbf{{{M}}}}[\mathbf{U},\hat{\mathbf{U}}]$   whenever  ${\mathbf{{{M}}}}[\mathbf{U},\hat{\mathbf{U}}]\in \mathbf{Y}_n, \, \forall\,\mathbf{U},\hat{\mathbf{U}} \in \mathbf{Y}_n.$ Thus,
% we can reformulate the invariant condition \eqref{invariant condition} as
if $$
u_k=\sum\limits_{i=1}^{n_k}\delta_{k,i}\psi_{k,i}(x) \text{ and } \hat{u}_k=\sum\limits_{i=1}^{n_k}\hat{\delta}_{k,i}\psi_{k,i}(x),$$
then  there exist arbitrary functions  $\Lambda_{k,i},$ and linear functions  $\hat{\Lambda}_{k,i},$ for $i=1,2,\dots,n_k,$  such that
\begin{eqnarray}\label{invariant condition-linear}
	\begin{aligned}
			{\mathcal{{{M}}}_k}\left[ {u}_k,\hat{u}_k \right]
	 =&
		\sum\limits_{i=1}^{n_k}
		\Lambda_{k,i}\Big(
		\delta_{1,1},
		%\hat{\delta}_{1,1}(t-\tau_1),
		\dots, \delta_{1,n_1}
		\delta_{2,1},
		%	\hat{\delta}_{2,1}(t-\tau_2),
		\dots
		,%\hat{\delta}_{1,n_1}(t-\tau_1),
		\delta_{2,n_2}
		%\hat{\delta}_{2,n_2}(t-\tau_2)
		\Big)\psi_{k,i}(x)
		\\&	+
		\sum\limits_{i=1}^{n_k}
		\hat{	\Lambda}_{k,i}\Big(
		\eta_{k,1}\hat{\delta}_{1,1}, \eta_{k,2}	\hat{\delta}_{2,1},
		\dots, \eta_{k,1}\hat{\delta}_{1,n_1},
			\eta_{k,2}\hat{\delta}_{2,n_2}
		\Big)\psi_{k,i}(x)
		,k=1,2.
	\end{aligned}
\end{eqnarray}
The following  theorem explains the separable form of analytical solutions for the given system \eqref{1+1:delay} with time delays, in which time delay occurs only in linear terms of $\hat{u}_1 $ and $\hat{u}_2. $
\begin{theorem}
	{Assume that the {coupled nonlinear time-fractional system}  of PDEs  with  linear time delays is} of the form
	\begin{eqnarray}\label{linear 1+1 delay}
		\left( 	{\partial^{\alpha_1}_{t}} u_1,
		{\partial^{\alpha_2}_{t}} u_2
		\right)
		={\mathbf{{{M}}}}[\mathbf{U},\hat{\mathbf{U}}],\alpha_1,\alpha_2>0,
	\end{eqnarray}
	where  $	{\partial^{\alpha_k}_{t}}(\cdot)$ represents  the fractional partial  time derivatives of order $\alpha_k,k=1,2,$ in the sense of either Caputo ${^C\partial^{\alpha_k}_{t}}(\cdot)$ or Riemann-Liouville ${^{RL}\partial^{\alpha_k}_{t}}(\cdot)$  as given in \eqref{Caputo}and \eqref{RL}, respectively, and
	${\mathbf{{{M}}}}[\mathbf{U},\hat{\mathbf{U}}]$ is given  in \eqref{linear-delay operator}. Then, the given {coupled nonlinear time-fractional system} \eqref{linear 1+1 delay} admits the following generalized separable analytical solution
		$$\begin{aligned}\mathbf{U}=(u_1,u_2),&u_k:=u_k(x,t)=\sum\limits_{i=1}^{n_k}\delta_{k,i}(t)\psi_{k,i}(x),
	\text{along with }
	\\	D^{\alpha_k}_t{\delta_{k,i}(t)}=&
		\Lambda_{k,i}\Big(
		\delta_{1,1}(t),
		%\hat{\delta}_{1,1}(t-\tau_1),
		\delta_{2,1}(t),
		%	\hat{\delta}_{2,1}(t-\tau_2),
		\dots, \delta_{1,n_1}(t)
		,%\hat{\delta}_{1,n_1}(t-\tau_1),
		\delta_{2,n_2}(t)
		%\hat{\delta}_{2,n_2}(t-\tau_2)
		\Big)
		 +
		\sum\limits_{i=1}^{n_k}
		\hat{	\Lambda}_{k,i}\Big(
		\eta_{k,1}
		\hat{\delta}_{1,1}(t-\tau_1),\\&		\eta_{k,2}	\hat{\delta}_{2,1}(t-\tau_2),
		 	\ldots, \eta_{k,1}\hat{\delta}_{1,n_1}(t-\tau_1),
		\eta_{k,2}\hat{\delta}_{2,n_2}(t-\tau_2)
		\Big),k=1,2,
	\end{aligned}
	$$ whenever $\mathbf{Y}_n$ is invariant vector space of  ${\mathbf{{{M}}}}[\mathbf{U},\hat{\mathbf{U}}]$ given in \eqref{linear-delay operator}.
	\label{corollary:thm1}
\end{theorem}
\begin{proof}
	Similar to the proof of Theorem \ref{main result}.
\end{proof}
Next, we present the applicability of the invariant subspace method through a {{coupled nonlinear time-fractional DR system}} with time delays in the following section.
\section{Invariant vector spaces  and solutions of the given system \eqref{DREs}}\label{sec3}
Consider a {{coupled nonlinear time-fractional DR system}} with multiple linear time delays:
\begin{eqnarray}
		{\partial^{\alpha_k}_{t}} u_k
	=&\mathcal{R}_k(\mathbf{U},\hat{\mathbf{U}}),0<\alpha_1,\alpha_2\leq2,k=1,2,
\end{eqnarray}
where $\mathcal{R}_k(\mathbf{U},\hat{\mathbf{U}})$ is given in \eqref{components of dre}.
Here, we consider the {fractional partial time derivatives } ${\partial^{\alpha_k}_{t}}(\cdot)$ of order $\alpha_k$ in the sense of {either} Caputo ${^C\partial^{\alpha_k}_{t}}(\cdot)$ {or} Riemann-Liouville ${^{RL}\partial^{\alpha_k}_{t}}(\cdot)$  of $u_k,k=1,2.$
 \subsection{Invariant vector spaces $\mathbf{Y}_4$ for the given system \eqref{DREs} }
 This subsection explains how to find invariant vector spaces  {$\mathbf{Y}_n=Y_{1,n_1}\times Y_{2,n_2},$ } with $n=4$ and $n_1=n_2=2,$  for the given system \eqref{DREs}.
% Now, we demonstrate how to find  $\mathbf{Y}_4$   for $\mathbf{R}[\mathbf{U},\hat{\mathbf{U}}].$
 For this case,  $ H_k[\cdot]=0$ given in \eqref{lhodes} can be read as follows:
 \begin{eqnarray}
 	\begin{aligned}\label{lhodes:n_k=2}
 		H_k[\lambda(x)]	:=({\mathbf{D}^{2}_x}+\mu_{k,1}{\mathbf{D}_x}+\mu_{k,0})\lambda(x)=0,
 	\end{aligned}
 \end{eqnarray}
for $k=1,2.$ It is known that $\mathbf{Y}_4$ is an invariant vector space for $\mathbf{R}[\mathbf{U},\hat{\mathbf{U}}]$ if  $\mathbf{R}[\mathbf{U},\hat{\mathbf{U}}]\in \mathbf{Y}_4= Y_{1,2}\times Y_{2,2}$ for all  $u_k\in Y_{k,2}=\mathcal{L}\{\psi_{k,1}(x),\psi_{k,2}(x)\},k=1,2.$
Thus, the invariant conditions  can be read as follows:
 \begin{eqnarray}
\begin{aligned}\label{lhodes in con n1=2}
	&
	\left( {\mathbf{D}^{2}_x}+\mu_{k,1}{\mathbf{D}_x}+\mu_{k,0}\right)\mathcal{R}_k[\mathbf{U},\hat{\mathbf{U}}]=0,\text{whenever}
  \left( {\mathbf{D}^{2}_x}+\mu_{k,1}{\mathbf{D}_x}+\mu_{k,0}\right) u_k =0,k=1,2.
\end{aligned}
 \end{eqnarray}
From the above-invariant condition \eqref{lhodes in con n1=2}, we get an over-determined system of equations.
	Solving the over-determined system of equations, we obtain various types of the four-dimensional invariant vector spaces $\mathbf{Y}_4$ and their corresponding   vector nonlinear delay differential operators $\mathbf{R}[\mathbf{U},\hat{\mathbf{U}}].$		Additionally, we note that the basis of the component vector spaces $ Y_{1,2}$ and $ Y_{2,2}$ can be determined in the form of   exponential, trigonometric, and polynomial functions with respect to  $x,$  which are discussed below in detail.
	\begin{case}\label{exp-case}
\textbf{Exponential invariant  vector spaces $\mathbf{Y}_4$ and their corresponding $\mathbf{R}[\mathbf{U},\hat{\mathbf{U}}]$:	}
The obtained exponential invariant vector spaces $\mathbf{Y}_4$ with their corresponding $\mathbf{R}[\mathbf{U},\hat{\mathbf{U}}]$ are given below.
	\begin{enumerate}
		\item \label{case1}
	$\mathbf{R}[\mathbf{U},\hat{\mathbf{U}}]$ admits $\mathbf{Y_4}=\mathcal{L}\{1,e^{\mu_{2,1}x}\}	\times\mathcal{L}\{1,e^{-\mu_{2,1}x}\}$	if
		$$\begin{array}{ll}
		 	\mathcal{R}_1(\mathbf{U},\hat{\mathbf{U}})
			=&\partial_x
		\big[
		(a_{12}u_2+a_{11}u_1+a_{10})
		\partial_x u_1
 +(b_{11}u_1+b_{10})
		\partial_x u_2
		\big ]
	-\mu_{{2,1}}^{2}(b_{{11}}u_1u_2
			\\&		+2a_{{11}}u_1^2
	+b_{{10}}u_2 )
			+\gamma_{11}u_1
		+\gamma_{10}
	 +\hat{\gamma}_{11}\hat{u}_1
	,
	\\
	\mathcal{R}_2(\mathbf{U},\hat{\mathbf{U}})= &\partial_x
	\big[
	(a_{22}u_2+a_{20})
	\partial_x u_1
	+(b_{22}u_2+b_{21}u_1+b_{20})
	\partial_x u_2
	\big]
		-\mu_{{2,1}}^{2}(a_{{22}}u_1u_2
\\&	+2b_{{22}}u_2^2
	+a_{20})u_1
		+
	\gamma_{22}u_2
	+\gamma_{20}
	+\hat{\gamma}_{22}\hat{u}_2.
\end{array}
	$$
\item \label{case2}
		$\mathbf{R}[\mathbf{U},\hat{\mathbf{U}}]$ admits $\mathbf{Y_4}=\mathcal{L}\{1,e^{-\mu_{1,1}x}\}	\times\mathcal{L}\{1,e^{-{\frac{
					\mu_{1,1} ( a_{{22}}+b_{{21}})x }{b_{{21}}}}}\}$	if
	$$\begin{array}{ll}
	\mathcal{R}_1(\mathbf{U},\hat{\mathbf{U}})= &\partial_x
		\big[
		(a_{12}u_2+a_{11}u_1+a_{10})
		\partial_x u_1
	 	+(b_{10}-\frac{a_{12}a_{22}}{a_{22}+b_{21}}u_1)
		\partial_x u_2
		\big]-2\mu_{{1,1}}^{2}a_{{11}}u_1^2
			\\&	-\frac{\mu_{{1,1}}^{2}(a_{22}+b_{21})^2b_{10}}{b_{21}^2}u_1u_2	+\gamma_{11}u_1
		+\gamma_{10}
	 	+\frac{\mu_{{1,1}}^{2}(a_{22}+b_{21})a_{12}a_{22}}{b_{21}^2}u_2 	+\hat{\gamma}_{11}\hat{u}_1
				,
	%	\end{aligned}
%$$
\\ %$$\begin{array}
	\mathcal{R}_2(\mathbf{U},\hat{\mathbf{U}})=&\gamma_{20}+\partial_x
		\big[
		(a_{22}u_2+a_{20})
		\partial_x u_1
	 +(b_{22}u_2+b_{21}u_1+b_{20})
		\partial_x u_2
		\big] 	
		+\gamma_{22}u_2
			\\& 	-\mu_{{1,1}}^{2}(a_{{22}}u_2
				+a_{20})u_1
 -\frac{2\mu_{{1,1}}^{2}(a_{22}+b_{21})^2b_{22}}{b_{21}^2}u_2^2
		+\hat{\gamma}_{22}\hat{u}_2.
	\end{array}
	$$
\item \label{case5}
$\mathbf{R}[\mathbf{U},\hat{\mathbf{U}}]$ admits $\mathbf{Y_4}=\mathcal{L}\{1,e^{-\mu_{1,1}x}\}	\times\mathcal{L}\{1,e^{-\mu_{2,1}x}\}$	if
$$\begin{array}{ll}
 	\mathcal{R}_1(\mathbf{U},\hat{\mathbf{U}})=& \partial_x
	\big[
	(a_{12}u_2+a_{11}u_1+a_{10})
	\partial_x u_1
 +(b_{10}-
	\frac{a_{12}(\mu_{1,1}+\mu_{2,1})}{\mu_{2,1}}u_1)
	\partial_x u_2
	\big]
\\& 	+\gamma_{11}u_1 +	\mu_{2,1}a_{12}(\mu_{2,1}+\mu_{1,1})u_1u_2
-2\mu_{1,1}^2a_{11}u_1^2	\\&-\mu_{2,1}^2b_{10}u_2
 +\hat{\gamma}_{11}\hat{u}_1+\gamma_{10}
	,
\\
	\mathcal{R}_2(\mathbf{U},\hat{\mathbf{U}})=&
a_{20}
\partial_{xx} u_1+	\partial_x\left[
	(b_{22}u_2+b_{20})
	\partial_x u_2
	\right]
	 +\gamma_{22}u_2	-2\mu_{{2,1}}^{2}{b_{22}}u_2^2
-\mu_{{1,1}}^{2}a_{20}u_1
\\&	+\gamma_{20}	+\hat{\gamma}_{22}\hat{u}_2.
\end{array}
$$
\item \label{case7}
$\mathbf{R}[\mathbf{U},\hat{\mathbf{U}}]$ admits $\mathbf{Y_4}=\mathcal{L}\{1,e^{-\frac{\mu_{1,1}}{2}x}\}	\times\mathcal{L}\{1,e^{-\mu_{2,1}x}\}$	if
$$\begin{array}{ll}
	\mathcal{R}_1(\mathbf{U},\hat{\mathbf{U}})=& \partial_x
	\big[
	(a_{11}u_1-\frac{2b_{11}}{3}u_2+a_{10})
	\partial_x u_1
 	+(b_{10}+b_{11}u_1)
	\partial_x u_2
	\big]	+\gamma_{11}u_1	
	\\&+	\mu_{2,1}^2(b_{11}u_1u_2-\dfrac{a_{11}}{2}u_1^2-b_{10}u_2)
+\gamma_{10}
+\hat{\gamma}_{11}\hat{u}_1	,
\\
 	\mathcal{R}_2(\mathbf{U},\hat{\mathbf{U}})=&
\partial_x
\big[
(a_{21}u_1-{3}b_{21}u_2
-\frac{4\gamma_{21}}{\mu_{2,1}^2})
\partial_x u_1
 +
	(b_{22}u_2+b_{20})
	\partial_x u_2
	\big]+\gamma_{22}u_2 \\
	&	+\gamma_{21}u_1
	+\gamma_{20}
	-\mu_{{2,1}}^{2}(\frac{3b_{21}}{4}u_1u_2+2{b_{22}}u_2^2+\frac{a_{21}}{8}u_1^2	)
	+\hat{\gamma}_{22}\hat{u}_2.
\end{array}
$$
\item \label{case11}
$\mathbf{R}[\mathbf{U},\hat{\mathbf{U}}]$ admits $\mathbf{Y_4}=\mathcal{L}\{1,e^{-{\mu_{1,1}}x}\}	\times\mathcal{L}\{e^{\mu_{1,1}x},e^{-\mu_{1,1}x}\}$	if
$$\begin{aligned}
	\mathcal{R}_1(\mathbf{U},\hat{\mathbf{U}})=& \partial_x
	\big[
	(a_{12}u_2+a_{11}u_1+a_{10})
	\partial_x u_1
 	+(b_{10}-2a_{12}u_1)
	\partial_x u_2
	\big]	+\gamma_{11}u_1
	+\hat{\gamma}_{11}\hat{u}_1
	\\&+	\mu_{1,1}^2(2a_{12}u_1u_2-2{a_{11}}u_1^2-b_{10}u_2)
+\gamma_{10}
	,
	\\
 	\mathcal{R}_2(\mathbf{U},\hat{\mathbf{U}})=
&	\partial_x
	\big[
	(a_{20}-b_{21}u_2
)
	\partial_x u_1
	+
	(b_{21}u_1
  +b_{20})
	\partial_x u_2
	\big]+\gamma_{22}u_2
			+\hat{\gamma}_{22}\hat{u}_2.
\end{aligned}
$$
\item \label{case13}
$\mathbf{R}[\mathbf{U},\hat{\mathbf{U}}]$ admits $\mathbf{Y_4}=\mathcal{L}\{1,e^{-{\mu_{1,1}}x}\}	\times\mathcal{L}\{e^{\frac{\mu_{1,1}}{2}x},{e^{\frac{\mu_{1,1}}{2}x}}\}$	if
$$\begin{array}{lll}
 	\mathcal{R}_1(\mathbf{U},\hat{\mathbf{U}})=&\partial_x
	\big[
	(a_{11}u_1+a_{10})
	\partial_x u_1
		\big]+b_{10}\partial_{xx} u_{2}
		 +\gamma_{11}u_1
		+\hat{\gamma}_{11}\hat{u}_1
 	\\&-	\mu_{1,1}^2(2{a_{11}}u_1^2+ \frac{b_{10}}{4}u_2)
	+\gamma_{10}
	,
	\\
	 \mathcal{R}_2(\mathbf{U},\hat{\mathbf{U}})= &\gamma_{22}u_2
	+
	\partial_x
	\big[
	(b_{22}u_2
	)
	\partial_x u_1
	+
	(b_{21}u_1   +b_{20})
	\partial_x u_2
	\big] -\frac{3\mu_{1,1}^2(2a_{22}+b_{21})}{4}u_1u_2+\hat{\gamma}_{22}\hat{u}_2.
\end{array}
$$	
\item \label{case15}
$\mathbf{R}[\mathbf{U},\hat{\mathbf{U}}]$ admits $\mathbf{Y_4}=\mathcal{L}\{1,e^{-{\mu_{1,1}}x}\}	\times\mathcal{L}\{e^{\mu_{1}x},{e^{\mu_{2}x}}\},$ $\mu_k=-\frac{\mu_{2,2}+(-1)^k\sqrt{\mu_{2,1}^2-4\mu_{2,0}}}{2},$ $k=1,2,$ whenever
$$\begin{array}{lll}
	\mathcal{R}_1(\mathbf{U},\hat{\mathbf{U}})=&\partial_x
	\big[
	(a_{11}u_1+a_{10})
	\partial_x u_1
	\big]
 +\gamma_{11}u_1	+\hat{\gamma}_{11}\hat{u}_1
	-	2\mu_{1,1}^2a_{11}u_1^2
	+\gamma_{10}
	,
\\
	\mathcal{R}_2(\mathbf{U},\hat{\mathbf{U}})=&
\gamma_{22}u_2
+	\partial_x
	\big[	(b_{21}u_1+b_{20})
	\partial_x u_2
 -	\frac{b_{21}(\mu_{1,1}+\mu_{2,1})}{\mu_{1,1}}u_2
\partial_x u_1
	\big]
	\\&+b_{21}(\mu_{1,1}^2 +\mu_{1,1}\mu_{2,1}
	+\mu_{2,0})u_1u_2+\hat{\gamma}_{22}\hat{u}_2.
\end{array}
$$
\item \label{case17}
$\mathbf{R}[\mathbf{U},\hat{\mathbf{U}}]$ admits $\mathbf{Y_4}=\mathcal{L}\{1,e^{-{\mu_{1,1}}x}\}	\times\mathcal{L}\{e^{\mu_{1,1}x},x{e^{\mu_{1,1}x}}\}$
	if
$$\begin{aligned}
&	\mathcal{R}_1(\mathbf{U},\hat{\mathbf{U}})=\partial_x
	\big[
	(a_{12}u_2+a_{11}u_1+a_{10})
	\partial_x u_1
	\big] +\gamma_{11}u_1
	+\hat{\gamma}_{11}\hat{u}_1
	-	2\mu_{1,1}^2a_{11}u_1^2
	+\gamma_{10}
	,
\\
&	\mathcal{R}_2(\mathbf{U},\hat{\mathbf{U}})=
	\partial_x
	\big[
{b_{21}}u_2
	\partial_x u_1
	+
	(b_{21}u_1   +b_{20})
	\partial_x u_2
	\big]+\gamma_{22}u_2
	+\hat{\gamma}_{22}\hat{u}_2.
\end{aligned}
$$
	
\item \label{case24}
$\mathbf{R}[\mathbf{U},\hat{\mathbf{U}}]$ admits $\mathbf{Y_4}=\mathcal{L}\{1,e^{-{\mu_{2,1}}x}\}	\times\mathcal{L}\{1,e^{-{\mu_{2,1}}x}\}$
if

$$\begin{aligned}
&	\mathcal{R}_k(\mathbf{U},\hat{\mathbf{U}})=\partial_x
	\big[
	(a_{k2}u_2+a_{k1}u_1+a_{k0})
	\partial_x u_1
	+	(b_{k2}u_2+b_{k1}u_1+b_{k0})
	\partial_x u_2
	\big]+ \gamma_{k1}u_1
	\\ &+\gamma_{k2}u_2	+\gamma_{k0}
	+\hat{\gamma}_{k1}\hat{u}_1	+\hat{\gamma}_{k2}\hat{u}_2
	-	2\mu_{1,1}^2[a_{k1}u_1^2   +b_{k2}u_2^2+2(a_{k2}+b_{k1})u_1u_2]
	,k=1,2.
\end{aligned}
$$	\end{enumerate}	
	
	\end{case}
\begin{case}\label{exp-pol case}
\textbf{Exponential-polynomial  invariant  vector spaces $\mathbf{Y}_4$ and their corresponding $\mathbf{R}[\mathbf{U},\hat{\mathbf{U}}]$:}
	The obtained exponential-polynomial invariant vector spaces $\mathbf{Y}_4$ with their corresponding $\mathbf{R}[\mathbf{U},\hat{\mathbf{U}}]$ are given below.
%Basis functions of $ Y_{1,2}$ and $ Y_{2,2}$ are in terms of exponential or polynomial functions in $x$:
	\begin{enumerate}
		\item \label{case3}
		$\mathbf{R}[\mathbf{U},\hat{\mathbf{U}}]$ admits $\mathbf{Y_4}=\mathcal{L}\{1,e^{\mu_{1,1}x}\}	\times\mathcal{L}\{1,x\}$	if
		$$\begin{array}{ll}
			\mathcal{R}_1(\mathbf{U},\hat{\mathbf{U}})=& \partial_x
			\big[
			(a_{11}u_1+a_{10})
			\partial_x u_1
			+(b_{11}u_1
			 +b_{10})
			\partial_x u_2
			\big]
				-2\mu_{{1,1}}^{2}a_{{11}}u_1^2
		 			+\gamma_{11}u_1
		\\&	+\gamma_{10}
			+\hat{\gamma}_{11}\hat{u}_1
			,
			\\
		\mathcal{R}_2(\mathbf{U},\hat{\mathbf{U}})= &\partial_x
			\big[
			(a_{20}-b_{21}u_2)
			\partial_x u_1
			+(b_{22}u_2 +b_{21}u_1+b_{20})
			\partial_x u_2
			\big]
				+\mu_{{1,1}}^{2}(b_{{21}}u_2
						-a_{20})u_1
			\\&+	\gamma_{22}u_2
		+\gamma_{20}
			+\hat{\gamma}_{22}\hat{u}_2.
		\end{array}
		$$
\item \label{case4}
		$\mathbf{R}[\mathbf{U},\hat{\mathbf{U}}]$ admits $\mathbf{Y_4}=	\mathcal{L}\{1,x\}\times\mathcal{L}\{1,e^{\mu_{2,1}x}\}$	if
				$$\begin{array}{ll}
		\mathcal{R}_1(\mathbf{U},\hat{\mathbf{U}})= &\gamma_{10}
		+ \partial_x
			\big[
			(a_{12}u_2+a_{11}u_1	  +a_{10})
			\partial_x u_1
		-(a_{12}u_1-b_{10})
			\partial_x u_2
			\big]
			\\&	+\mu_{{2,1}}^{2}(a_{{12}}u_1^2-b_{10}u_2)
			+\gamma_{11}u_1
			+\hat{\gamma}_{11}\hat{u}_1
			,
			\\
			\mathcal{R}_2(\mathbf{U},\hat{\mathbf{U}})= &\partial_x
			\big[
			(a_{20}+a_{22}u_2)
			\partial_x u_1
		 +(b_{22}u_2 +b_{20})
			\partial_x u_2
			\big]
			-2\mu_{{2,1}}^{2}b_{22}u_2^2+
			\gamma_{22}u_2
		\\&	+\gamma_{20}
		+\hat{\gamma}_{22}\hat{u}_2.
		\end{array}
		$$
		\item \label{case9}
		$\mathbf{R}[\mathbf{U},\hat{\mathbf{U}}]$ admits $\mathbf{Y_4}=	\mathcal{L}\{1,x\}\times\mathcal{L}\{e^{\mu_1 x},e^{\mu_2 x}\},\mu_k=-\frac{\mu_{2,2}+(-1)^k\sqrt{\mu_{2,1}^2-4\mu_{2,0}}}{2},$ $k=1,2,$	whenever
		$$\begin{array}{ll}
		\mathcal{R}_1(\mathbf{U},\hat{\mathbf{U}})=&\partial_x
			\left[
			(a_{11}u_1+a_{10})
			\partial_x u_1
					\right]
					 	+\gamma_{11}u_1
			+\gamma_{10}
			+\hat{\gamma}_{11}\hat{u}_1
			,
	\\ 	\mathcal{R}_2(\mathbf{U},\hat{\mathbf{U}})=&\partial_x
			\left[
			(a_{20}+a_{22}u_2)
			\partial_x u_1
						\right]   +b_{20}
						\partial_{xx} u_2
				+
			\gamma_{22}u_2
				+\hat{\gamma}_{22}\hat{u}_2.
		\end{array}
		$$
	\item \label{case10}
$\mathbf{R}[\mathbf{U},\hat{\mathbf{U}}]$ admits $\mathbf{Y_4}=	\mathcal{L}\{1,x\}\times\mathcal{L}\{e^{-\frac{\mu_{2,1}}{2} x},xe^{-\frac{\mu_{2,1}}{2} x}\}$	if
$$\begin{array}{ll}
	\mathcal{R}_1(\mathbf{U},\hat{\mathbf{U}})=&\partial_x
	[
	(a_{11}u_1+a_{10})
	\partial_x u_1
	]
	+\gamma_{11}u_1
	 +\hat{\gamma}_{11}\hat{u}_1+\gamma_{10}
	,
\\
\mathcal{R}_2(\mathbf{U},\hat{\mathbf{U}})=&\partial_x
	[
	(a_{20}+a_{22}u_2)
	\partial_x u_1
	+(b_{20} +b_{21}u_1)
	\partial_x u_2
	]
	-\frac{\mu_{2,1}^2 b_{21}}{4} u_1u_2
	+\gamma_{22}u_2
	+\hat{\gamma}_{22}\hat{u}_2.
\end{array}
$$	
	\item \label{case31}
$\mathbf{R}[\mathbf{U},\hat{\mathbf{U}}]$ admits $\mathbf{Y_4}=	\mathcal{L}\{e^{\mu_1 x},e^{\mu_2 x}\}\times\mathcal{L}\{1,x\},\mu_k=-\frac{\mu_{1,2}+(-1)^k\sqrt{\mu_{1,1}^2-4\mu_{1,0}}}{2},$ $k=1,2,$	whenever
$$\begin{array}{ll}
& \mathcal{R}_1(\mathbf{U},\hat{\mathbf{U}})=a_{10}\partial_{xx}u_1+\partial_x
	\left[
	(b_{10}+b_{11}u_1)
	\partial_x u_2
	\right]	 +
	\gamma_{11}u_1
	+\hat{\gamma}_{11}\hat{u}_1,
	\\
	&	\mathcal{R}_2(\mathbf{U},\hat{\mathbf{U}})=\partial_x
	\left[
	(b_{22}u_2+b_{20})
	\partial_x u_2
	\right]
	+\gamma_{22}u_2
 	+\gamma_{20}
	+\hat{\gamma}_{22}\hat{u}_2.
\end{array}
$$	
	\end{enumerate}
	
\end{case}
	\begin{case} \label{trigom-exp case}
\textbf{Exponential-trigonometric invariant  vector spaces $\mathbf{Y}_4$ and their corresponding $\mathbf{R}[\mathbf{U},\hat{\mathbf{U}}]:$	}
The obtained exponential-trigonometric invariant vector spaces $\mathbf{Y}_4$ with their corresponding $\mathbf{R}[\mathbf{U},\hat{\mathbf{U}}]$ are given below.
%Basis functions of	$ Y_{1,2}$ and $ Y_{2,2}$ are in terms of exponential or trigonometric functions in $x$:
\begin{enumerate}
 \item \label{case12}
	$\mathbf{R}[\mathbf{U},\hat{\mathbf{U}}]$ admits $\mathbf{Y_4}=\mathcal{L}\{1,e^{-{\mu_{1,1}}x}\}	\times\mathcal{L}\{sin(\sqrt{\mu_{2,0}}x),cos(\sqrt{\mu_{2,0}}x)\},$ $\mu_{2,0}>0,$ if
	
	$$\begin{array}{ll}
		\mathcal{R}_1(\mathbf{U},\hat{\mathbf{U}})= &\partial_x
		\big[
		(a_{11}u_1+a_{10})
		\partial_x u_1
			\big]+b_{10}
			\partial_{xx} u_2
		  -2\mu_{1,1}^2a_{11}u_1^2+	\mu_{2,0}b_{10}u_2
		+\gamma_{11}u_1
	\\&	+\hat{\gamma}_{11}\hat{u}_1+\gamma_{10}
		,
		\\
	 	\mathcal{R}_2(\mathbf{U},\hat{\mathbf{U}})=&
	\gamma_{22}u_2
	+	\partial_x
		\big[
				(b_{21}u_1+b_{20})
				\partial_x u_2
			 	-(b_{21}u_2)
			\partial_x u_1		\big]	+b_{21}(\mu_{{1,1}}^2+\mu_{{2,0}})u_1u_2
		\\&	+\hat{\gamma}_{22}\hat{u}_2.
	\end{array}
	$$
\item \label{case29}
	$\mathbf{R}[\mathbf{U},\hat{\mathbf{U}}]$ admits $\mathbf{Y_4}=\mathcal{L}\{sin(\sqrt{\mu_{1,0}}x),cos(\sqrt{\mu_{1,0}}x)\}	\times\mathcal{L}\{1,e^{-{\mu_{2,1}}x}\},$ $\mu_{1,0}>0,$ whenever
	$$\begin{array}{ll}
 	\mathcal{R}_1(\mathbf{U},\hat{\mathbf{U}})=&
\hat{\gamma}_{22}\hat{u}_2	+	\partial_x
	\big[
	(a_{10}-b_{11}u_2)
	\partial_x u_1  	+
	(b_{11}u_1)
	\partial_x u_2
	\big]	+\gamma_{22}u_2
\\&	-b_{11}(\mu_{{2,1}}^2+\mu_{{1,0}})u_1u_2
			,
\\
	 \mathcal{R}_2(\mathbf{U},\hat{\mathbf{U}})= &\partial_x
	\big[
	(a_{21}u_1+a_{20})
	\partial_x u_1 +(b_{22}u_2   +b_{20})\partial_{x} u_2
	\big]	
+	\gamma_{{22}}u_2
	\\& 	+2\mu_{1,0}a_{21}u_1^2-2\mu_{2,1}^2b_{22}u_2^2
+\mu_{{1,0}}a_{20}u_1		+\hat{\gamma}_{11}\hat{u}_2+\gamma_{20}	.
	\end{array}
	$$
\end{enumerate}		
		
	\end{case}
\begin{case} \label{tro-polcase}
	\textbf{Polynomial-trigonometric invariant  vector spaces $\mathbf{Y}_4$ and their corresponding $\mathbf{R}[\mathbf{U},\hat{\mathbf{U}}]:$}
	The obtained polynomial-trigonometric invariant vector spaces $\mathbf{Y}_4$ with their corresponding $\mathbf{R}[\mathbf{U},\hat{\mathbf{U}}]$ are given below.
%Basis functions of	$ Y_{1,2}$ and $ Y_{2,2}$ are in terms of polynomial or trigonometric functions in $x$:
	\begin{enumerate}
		\item \label{case19} 	$\mathbf{R}[\mathbf{U},\hat{\mathbf{U}}]$ admits $	\mathbf{Y}_4=\mathcal{L}\{1,x\}\times \mathcal{L}\left\{sin(\sqrt{\mu_{2,0}}x),cos(\sqrt{\mu_{{2,0}}}x)\right\},$ $\mu_{2,0}>0,$ whenever
$$	\begin{array}{ll}
 \mathcal{R}_1(\mathbf{U},\hat{\mathbf{U}})=&
		\partial_x\Big[
	(a_{12}u_2+a_{11}u_1+a_{10})\partial_xu_1
 	+\left( b_{12}u_2-a_{12}u_1
	+\frac{\gamma_{{12}}}{\mu_{{2,0}}}\right) \partial_x u_2 \Big]+\gamma_{10}
	\\
&	-\mu_{{2,0}}a_{12}u_1u_2
+2\mu_{{2,0}}b_{12}u_2^2+\gamma_{12}u_2  +\gamma_{11}u_1	+
	\hat{\gamma}_{11}\hat{u}_1,
\\
 \mathcal{R}_2(\mathbf{U},\hat{\mathbf{U}})=&
	\partial_x
	\big[
	(a_{22}u_2+a_{20})\partial_x u_1+(b_{21}u_1 	  +b_{20})\partial_x u_2 \big]
	+\mu_{{2,0}}b_{21}u_1u_2
	+\gamma_{22}u_2
	+\hat{\gamma}_{22}\hat{u}_2,
\end{array}$$
\item\label{case27}
$\mathbf{R}[\mathbf{U},\hat{\mathbf{U}}]$ admits $	\mathbf{Y}_4= \mathcal{L} \{sin(\sqrt{\mu_{1,0}}x),$ $cos(\sqrt{\mu_{{1,0}}}x) \} $ $\times \mathcal{L}\{1,x\},$ $\mu_{1,0}>0,$
whenever
$$
\begin{aligned}
 	\mathcal{R}_1(\mathbf{U},\hat{\mathbf{U}})=&
	\partial_x
	\big[
	(a_{12}u_2+a_{10})\partial_x u_1+(b_{11}u_1
 +b_{10})\partial_x u_2 \big]
	+\mu_{{1,0}}a_{12}u_1u_2
			+\gamma_{11}u_1
	+\hat{\gamma}_{11}\hat{u}_1,
	\\
	\mathcal{R}_2(\mathbf{U},\hat{\mathbf{U}})=&
		\partial_x[
	(a_{21}u_1-b_{21}u_2+a_{20})\partial_xu_1
 	+( b_{22}u_2+b_{12}u_1
	+b_{20}) \partial_x u_2 ]+\gamma_{12}u_2
	+\gamma_{20}
	\\&	+\mu_{{1,0}}(2a_{21}u_1-b_{21}u_2
	+a_{20})u_1+
	\hat{\gamma}_{22}\hat{u}_2.
\end{aligned}$$
	\end{enumerate}
	
\end{case}
\begin{case} \label{poly case}
\textbf{Polynomial invariant  vector spaces $\mathbf{Y}_4$ and their corresponding $\mathbf{R}[\mathbf{U},\hat{\mathbf{U}}]$:}
	The obtained polynomial invariant vector space $\mathbf{Y}_4$ with their corresponding $\mathbf{R}[\mathbf{U},\hat{\mathbf{U}}]$ are given below.
%	Basis functions of	$ Y_{1,2}$ and $ Y_{2,2}$ are   polynomial functions in $x$:
	\begin{itemize}
		\item[$\bullet$]\label{case33}
$\mathbf{R}[\mathbf{U},\hat{\mathbf{U}}]$ admits 	$\mathbf{Y}_4= \mathcal{L}\{1,x\} \times \mathcal{L}\{1,x\}$ if
	$$	\begin{aligned}
	 	\mathcal{R}_k(\mathbf{U},\hat{\mathbf{U}})=&
		\partial_x\big[(a_{12}u_2+a_{11}u_1+a_{10})\partial_x u_1
		 +(b_{k2}u_2+b_{k1}u_1+b_{k0})\partial_x u_2 \big]+\gamma_{k2}u_2
		\\
		&	+\gamma_{k1}u_1+\gamma_{k0}
		+
		\hat{\gamma}_{k1}\hat{u}_1+\hat{\gamma}_{k2}\hat{u}_2,
		k=1,2.
	\end{aligned}
$$
\end{itemize}
\end{case}
%%%%%%
%%%%%
%%%%lambda_1=\mu_{{1,1}},c_{}=gamma_{{}} eta_{i}=\mu_{2i}
Next, we explain  the procedure to derive analytical solutions for the given system \eqref{DREs} using the above obtained
invariant vector spaces $\mathbf{Y}_4$ through the following examples.
\subsection{Solutions of the given system \eqref{DREs}}
\begin{example}
	\label{eg1}
	Consider a {coupled nonlinear time-fractional DR system} with linear time delays as follows:
	\begin{eqnarray}
		\begin{aligned}\label{DRE:eg1}
		{\partial^{\alpha_1}_{t}} u_1=&	%\mathcal{R}_1(\mathbf{U},\hat{\mathbf{U}})\equiv
				\partial_x\left[
		(a_{12}u_2+a_{11}u_1+a_{10})\partial_xu_1
 	+\left( b_{12}u_2-a_{12}u_1
		+\frac{\gamma_{{12}}}{\mu_{{2,0}}}\right) \partial_x u_2 \right]
			\\&	-\mu_{{2,0}}a_{12}u_1u_2
			+2\mu_{{2,0}}b_{12}u_2^2+\gamma_{12}u_2  +\gamma_{11}u_1	+
		\gamma_{10}
		+	\hat{\gamma}_{11}\hat{u}_1,	
		\\	{\partial^{\alpha_2}_{t}} u_2=&	%\mathcal{R}_k(\mathbf{U},\hat{\mathbf{U}})\equiv
		\partial_x
		\left[
		(a_{22}u_2+a_{20})\partial_x u_1+(b_{21}u_1+b_{20})\partial_x u_2 \right]
		+\mu_{{2,0}}b_{21}u_1u_2
			+\gamma_{22}u_2
		+\hat{\gamma}_{22}\hat{u}_2,
		\\
		\hat{u}_k(x,t)=&\phi_k(x,t),t\in[-\tau_k,0],\tau_k>0,k=1,2, 	
	\end{aligned}
	\end{eqnarray}
	where $ 0<\alpha_1,\alpha_2\leq2,$ ${\partial^{\alpha_k}_{t}}(\cdot) $ is  the  {fractional partial time derivative} of order $\alpha_k,k=1,2,$ in the sense of {either} Caputo ${^C\partial^{\alpha_k}_{t}}(\cdot)$ {or} Riemann-Liouville ${^{RL}\partial^{\alpha_k}_{t}}(\cdot)$   given in \eqref{Caputo} and \eqref{RL}, respectively.
	Here, we consider  the initial and the boundary conditions as given in \eqref{ic-bc caputo}  and \eqref{ic-bc rl} under ${^C\partial^{\alpha_k}_{t}}(\cdot)$ and ${^{RL}\partial^{\alpha_k}_{t}}(\cdot),$ respectively.
From the Case \ref{tro-polcase}.\ref{case19}  in the previous subsection,   the above-{system} \eqref{DRE:eg1} admits the invariant vector space
\begin{eqnarray}
	\label{vector space example 1}
	\mathbf{Y}_4=\mathcal{L}\{1,x\}\times \mathcal{L} \{sin(\sqrt{\mu_{2,0}}x),cos(\sqrt{\mu_{{2,0}}}x)\},
\end{eqnarray}
where $\mu_{2,0}>0.$ Using the above invariant vector space $	\mathbf{Y}_4$, we get  the analytical solution for  \eqref{DRE:eg1} as follows:
\begin{eqnarray}
	\begin{aligned}\label{form of sol eg1}
			u_1(x,t)=&\delta_{1,1}(t)+\delta_{1,2}(t)x, \text{ and  }
	u_2(x,t)=	\delta_{2,1}(t)sin(\sqrt{\mu_{2,0}}x)
	 	+\delta_{2,2}(t)cos(\sqrt{\mu_{2,0}}x),
		\end{aligned}
\end{eqnarray}
where the unknown functions $\delta_{k,i}(t),k,i=1,2$ are to be determined.
% by  a system of time-fractional ODEs based on the invariant vector space $\mathbf{Y}_4$ given in \eqref{vector space example 1}.
For $a_{11}=b_{12}=0$ and $a_{22}=-b_{21},$  the given system \eqref{DRE:eg1} reduces into the following system of fractional-order time derivative delay ODEs:
\begin{eqnarray}
	\begin{aligned}
	\label{eq1:ods1}
&{{D}^{\alpha_1}_{t}}\delta_{1,1}(t)=\gamma_{{10}}+\gamma_{{11}}\delta_{1,1}(t)+\hat{\gamma}_{11}\hat{\delta}_{1,1}(t-\tau_1),
%	\label{eq1:ods2}
\, \,{{D}^{\alpha_1}_{t}}\delta_{1,2}(t)=\gamma_{{11}}\delta_{1,1}(t)+\hat{\gamma}_{11}\hat{\delta}_{1,1}(t-\tau_1), \\ &\&\
%	\label{eq1:ods3}
{{D}^{\alpha_2}_{t}}\delta_{2,i}(t)=
C\delta_{2,i}(t)+\hat{\gamma}_{22}\hat{\delta}_{2,i}(t-\tau_2),i=1,2,
\end{aligned}
\end{eqnarray}
where $C=\gamma_{{22}}-\mu_{20}.$ We solve the above-obtained  system of fractional-order time derivative delay ODEs \eqref{eq1:ods1} using   the Laplace transformation technique if ${{D}^{\alpha_k}_{t}},k=1,2$ is  fractional ordinary time derivative of order $\alpha_k$ in the sense of both the  Caputo ${^CD}^{\alpha_k}_t(\cdot)$ and the Riemann-Liouville ${^{RL}D}^{\alpha_k}_t(\cdot).$ So, we recall  the Laplace transformation of the  Caputo ${^CD}^{\alpha}_t(\cdot)$ and the Riemann-Liouville $ {^{RL}D}^{\alpha}_t(\cdot)$  fractional ordinary time derivatives of order $\alpha,\alpha\in(n-1,n],n\in\mathbb{N},$  as follows \cite{Diethelm2010}:
	$$
	\begin{array}{lll}
		&\mathbf{L}\left({{^CD}^{\alpha}_tf(t) };s\right)=& s^\alpha\mathbf{L}(f(t);s)  -\sum\limits_{q=0}^{n-1}s^{\alpha-(q+1)} {\mathbf{D}^q_t[f(t)]}\Big{|}_{t=0},
	\end{array}$$
		$$
		 \begin{array}{lll}
		  \text{and }&	\mathbf{L}\left({{^{RL}D}^{\alpha}_tf(t) };s\right)=& s^\alpha\mathbf{L}(f(t);s)  -\sum\limits_{q=0}^{n-1}s^{q}\,  {{^{RL}D}^{\alpha-q-1}_t[f(t)]}\Big{|}_{t=0},
	\end{array}$$
respectively, where $ \text{Re}(s)>0.$
Now, let us first consider the  fractional ordinary time derivative in the sense of Caputo, that is,  ${D^{\alpha_k}_{t}}(\cdot)={^CD}^{\alpha_k}_t(\cdot),0<\alpha_k\leq 2, k=1,2.$
Thus, we apply  the  Laplace transformation to the first equation of \eqref{eq1:ods1}.
{For }$0< \alpha_1,\alpha_2\leq1,$ we obtain
\begin{eqnarray}
&	\begin{aligned}\label{eq1alpha01}
&	s^{\alpha_1}\mathbf{L}(\delta_{1,1}(t);s)
	-\beta_{1,1}s^{\alpha_1-1}
	=\frac{\gamma_{{10}}}{s} \ +\gamma_{{11}}\mathbf{L}(\delta_{1,1}(t);s)
	 +\hat{\gamma}_{11}\mathbf{L}(\hat{\delta}_{1,1}(t-\tau_1);s).
	\end{aligned}
\end{eqnarray}
{For }$1< \alpha_1,\alpha_2\leq2,$ we obtain
\begin{eqnarray}
	&	\begin{aligned}\label{eq2alpha11}	&s^{\alpha_1}\mathbf{L}(\delta_{1,1}(t);s)
		-\beta_{1,1}s^{\alpha_1-1}	-\kappa_{1,1}s^{\alpha_1-2}
 	=\frac{\gamma_{{10}}}{s} +\gamma_{{11}}\mathbf{L}(\delta_{1,1}(t);s)   +\hat{\gamma}_{11}\mathbf{L}(\hat{\delta}_{1,1}(t-\tau_1);s),
			\end{aligned}
\end{eqnarray}
where
$\beta_{1,1}=\delta_{1,1}(0),$ $\kappa_{1,1}=
\mathbf{D}_t[\delta_{1,1}(t)]\Big{|}_{t=0} .$
After algebraic simplifications of  the above equations  \eqref{eq1alpha01}-\eqref{eq2alpha11} and %-\eqref{eq2alpha12}
 applying  the inverse Laplace transformation, we obtain the function $\delta_{1,1}(t) $ as follows:
\begin{eqnarray}
	\begin{aligned}
	\delta_{1,1}(t)=&\left\{
		\begin{array}{lll}
\hat{\gamma}_{11}F_{12}(t)\ast \varphi_{1,1}(t-\tau_1)\rho_{1}(t-\tau_1)
 +\gamma_{10}F_{11}(t)	
+	\beta_{1,1}F_{13}(t),
\\
\text{if } 0< \alpha_1\leq1,0<\alpha_2\leq2,
\\
\hat{\gamma}_{11}F_{12}(t)\ast \varphi_{1,1}(t-\tau_1)\rho_{1}(t-\tau_1)
 +\gamma_{10}F_{11}(t)	
+	\beta_{1,1}F_{13}(t)+	\kappa_{1,1}F_{14}(t),
\\ \text{if } 1< \alpha_1\leq2,0<\alpha_2\leq2.
	\end{array}
	\right.
	\label{delta11eg1}
		\end{aligned}
\end{eqnarray}
Similarly, we obtain
\begin{eqnarray}	\begin{aligned}	
\delta_{1,2}(t)=&\left\{
\begin{array}{lll}
\hat{\gamma}_{11}	F_{12}(t)\ast \varphi_{1,2}(t-\tau_1)\rho_{1}(t-\tau_1)
	  +	\beta_{1,2}F_{13}(t), \text{ if } 0< \alpha_1\leq1,
	0<\alpha_2\leq2,
	\\
\hat{\gamma}_{11}	F_{12}(t)\ast \varphi_{1,2}(t-\tau_1)\rho_{1}(t-\tau_1)
	 +	\beta_{1,2}F_{13}(t)+	\kappa_{1,2}F_{14}(t),
	\\
	\text{if } 1< \alpha_1\leq2,0<\alpha_2\leq2,
\end{array}
\right.
\label{delta12eg1}
%\end{aligned}
%\end{eqnarray}
%\begin{eqnarray}
%	\begin{aligned}	
	\\\delta_{2,i}(t)=&\left\{
	\begin{array}{lll}
			\hat{\gamma}_{22}F_{21}(t)\ast \varphi_{2,i}(t-\tau_2)\rho_{1}(t-\tau_2)
		 +	\beta_{2,i}F_{22}(t), \text{if } 0< \alpha_1\leq2,0<\alpha_2\leq1,
		\\
	\hat{\gamma}_{22}	F_{21}(t)\ast \varphi_{2,i}(t-\tau_2)\rho_{1}(t-\tau_2)
		 +\beta_{2,i}F_{22}(t)+	\kappa_{2,i}F_{23}(t), \\ \text{if } 0< \alpha_1\leq2,1<\alpha_2\leq2.
	\end{array}
	\right.	\label{delta2ieg1}
\end{aligned}
	\end{eqnarray}
Here
$\beta_{k,i}=\delta_{k,i}(0) , $ $\kappa_{k,i}=\mathbf{D}_t[\delta_{k,i}(t)]\Big{|}_{t=0},i=1,2,$   \begin{eqnarray}\label{f1i}
	\begin{aligned}
	&	F_{11}(t)=	\sum\limits_{m=0}^{\lfloor\frac{t}{\tau_{1}}\rfloor}\hat{\gamma}_{11}^m\eta_1^{\alpha_1(m+1)}  E^{m+1}_{\alpha_1,\alpha_1(m+1)+1}(\gamma_{{11}} \eta_1^{\alpha_1}), \delta_{k,i}(t)=\varphi_{k,i}(t),t\in[-\tau_{k},0],
\\&F_{12}(t)=	\sum\limits_{m=0}^{\lfloor\frac{t}{\tau_{1}}\rfloor}\hat{\gamma}_{11}^m\eta_1^{\alpha_1(m+1)-1}   E^{m+1}_{\alpha_1,\alpha_1(m+1)}(\gamma_{{11}} \eta_1^{\alpha_1}),\varphi_{k,i}(t)=\varphi_{k,i}(0),t>0,
		\\ &	F_{13}(t)=	\sum\limits_{m=0}^{\lfloor\frac{t}{\tau_{1}}\rfloor}\hat{\gamma}_{11}^m\eta_1^{\alpha_1m} E^{m+1}_{\alpha_1,\alpha_1m+1}(\gamma_{{11}}\eta_1^{\alpha_1}),
		 F_{14}(t)=	\sum\limits_{m=0}^{\lfloor\frac{t}{\tau_{1}}\rfloor}\hat{\gamma}_{11}^m\eta_1^{\alpha_1m+1}   E^{m+1}_{\alpha_1,\alpha_1m+2}(\gamma_{{11}}\eta_1^{\alpha_1}),
	\\
&	F_{22}(t)=	\sum\limits_{m=0}^{\lfloor\frac{t}{\tau_{2}}\rfloor}\hat{\gamma}_{22}^m\eta_2^{\alpha_2m}  E^{m+1}_{\alpha_2,\alpha_2m+1}\big(C\eta_2^{\alpha_2}\big),
	F_{23}(t)=	\sum\limits_{m=0}^{\lfloor\frac{t}{\tau_{2}}\rfloor}\hat{\gamma}_{22}^m\eta_2^{\alpha_2m+1}   E^{m+1}_{\alpha_2,\alpha_2m+2}\big(C\eta_2^{\alpha_2}\big),
\\	&F_{21}(t)= 	\sum\limits_{m=0}^{\lfloor\frac{t}{\tau_{2}}\rfloor}\hat{\gamma}_{22}^m\eta_2^{\alpha_2(m+1)-1}   E^{m+1}_{\alpha_2,\alpha_2(m+1)}\big(C\eta_2^{\alpha_2}\big),
\eta_k=t-\tau_{k}m, \rho_{1}(t)=\left\{\begin{array}{ll}
	0,&\text{if }t\geq0,
	\\
	1&\text{if }t<0,
\end{array} \right.&	
		\end{aligned}
	\end{eqnarray}
     $\ast$ {is the convolution product of functions,} {and} $  E^{\gamma_3}_{\gamma_1,\gamma_2}(t)=\sum\limits_{p=0}^\infty\dfrac{\Gamma(\gamma_3+p) t^{\gamma_1}}{p!\Gamma(\gamma_3)\Gamma(\gamma_1 p+\gamma_2)},$ $\gamma_3\in\mathbb{R}, \gamma_i>0,i,k=1,2,	 $ is the {generalized {3-parameter} Mittag-Leffler function} \cite{Diethelm2010}.
Thus, the obtained solution for the given {coupled nonlinear time-fractional DR system} \eqref{DRE:eg1} with the Caputo  {fractional partial time derivative} ${^C\partial^{\alpha_k}_t}(\cdot)$ as follows:
\begin{eqnarray*}
	\begin{aligned}
		&	u_1(x,t)=\delta_{1,1}(t)+\delta_{1,2}(t)x, \text{and }
			u_2(x,t)=\delta_{2,1}(t)sin(\sqrt{\mu_{2,0}}x)+\delta_{2,2}(t)cos(\sqrt{\mu_{2,0}}x),
	\end{aligned}
\end{eqnarray*}
where $\delta_{k,i}(t),i,k=1,2,$ are given in \eqref{delta11eg1}-\eqref{delta12eg1}.
%-\eqref{delta12eg1}.
The above-analytical solution of   \eqref{DRE:eg1} in the sense of the Caputo {fractional partial time derivative} ${^C\partial^{\alpha_k}_t}(\cdot)$ obeys the following initial conditions:
\begin{eqnarray*}
&\begin{aligned}
	&	u_1(x,0)=\beta_{1,1}+\beta_{1,2}x,
	u_2(x,0)=\beta_{2,1}sin(\sqrt{\mu_{2,0}}x)  +\beta_{2,2}cos(\sqrt{\mu_{2,0}}x), %\text{and }
  \\ &\partial_t[u_1(x,t)]\Big{|}_{t=0}=\kappa_{1,1}
+\kappa_{1,2}x, \,\&\,
%\end{aligned}
%\begin{aligned}
 \partial_t[u_2(x,t)]\Big{|}_{t=0}=\kappa_{2,1}sin(\sqrt{\mu_{2,0}}x)
 +\kappa_{2,2}cos(\sqrt{\mu_{2,0}}x),
\end{aligned}
`\\&
\begin{aligned}
	\text{and the boundary conditions: }		&u_1(0,t)=\delta_{1,1}(t) %\text{and }
		, 	 u_1(\lambda,t)=\delta_{1,1}(t)+\lambda\delta_{1,2}(t),	\\ 	&
u_2(0,t)=\delta_{2,2}(t),	\,\&\,
			u_2(\lambda,t)=\lambda_1\delta_{2,1}(t) +\lambda_2\delta_{2,2}(t),
		\end{aligned}
	\end{eqnarray*}
where $\delta_{k,i},i,k=1,2,$ are given in \eqref{delta11eg1},
%-\eqref{delta12eg1}
 $\lambda_1=sin(\lambda\sqrt{\mu_{2,0}})$ and $\lambda_2=cos(\lambda\sqrt{\mu_{2,0}}).$
%%%%%%%%%%%%%%%%%%%%%%%%%%%%%%%%%%%%%%%%%%%%%%%%%%%%%%%%%%%%%%%%%%%%%%%%%%%%%%%RL%%%%%%%%%%%%%%%%%%%%%%%%%%%%%%%%%%%%%%%%%%%%%%%%%%%%%%%%%%%%%%%%%%%%%%%%%%%%%%%%%%%%%%

Next, we consider the {fractional partial time derivative} of \eqref{DRE:eg1} in the sense of Riemann-Liouville. So, ${D^{\alpha_k}_{t}}(\cdot)={^{RL}{D}}^{\alpha_k}_{t}(\cdot),0<\alpha_k\leq 2, k=1,2,$  given in \eqref{eq1:ods1}.
We apply the Laplace transformation of the first equation in the obtained  equations \eqref{eq1:ods1}. Thus,
{for }$0< \alpha_1,\alpha_2\leq1,$ we obtain
	\begin{eqnarray}
	\begin{aligned}		\label{rl-eq1alpha01}
		%	(a)\, &
			s^{\alpha_1}\mathbf{L}(\delta_{1,1}(t);s)
			-\hat{\beta}_{1,1}
					=&\frac{\gamma_{{10}}}{s}
					+\gamma_{{11}}\mathbf{L}(\delta_{1,1}(t);s) 	  +\hat{\gamma}_{11}\mathbf{L}(\hat{\delta}_{1,1}(t-\tau_1);s) .
		\end{aligned}
	\end{eqnarray}
Also, {for }$1< \alpha_1,\alpha_2\leq2,$ we get
	\begin{eqnarray}
		\begin{aligned}\label{rl-eq1alpha12}%	(a)\,
	&	s^{\alpha_1}\mathbf{L}(\delta_{1,1}(t);s)
				-\hat{\beta}_{1,1}	-s 	\hat{\kappa}_{1,1}%{^{RL}D}^{\alpha_1-2}_{t}[\delta_{1,1}(t)]\Big{|}_{t=0}
 	=\gamma_{{11}}\mathbf{L}(\delta_{1,1}(t);s)+\frac{\gamma_{{10}}}{s}  +\hat{\gamma}_{11}\mathbf{L}(\hat{\delta}_{1,1}(t-\tau_1);s)
 .
		\end{aligned}
	\end{eqnarray}
 where $\hat{\beta}_{1,1}={^{RL}D}^{\alpha_1-1}_{t}[\delta_{1,1}(t)]\Big{|}_{t=0}$ and  $\hat{\kappa}_{1,1}={^{RL}D}^{\alpha_1-2}_{t}[\delta_{1,1}(t)]\Big{|}_{t=0}.$
After algebraic simplifications and applying the inverse Laplace transformation of the  equations \eqref{rl-eq1alpha01}-\eqref{rl-eq1alpha12}, we obtain $\delta_{1,1}(t)$ as follows:

\begin{eqnarray}
	\begin{aligned}
	\delta_{1,1}(t)=&\left\{
	\begin{array}{lll}
		\gamma_{10}{F}_{11}(t)	
		+  \varphi_{1,1}(t-\tau_1)\rho_{1}(t-\tau_1)
		\ast \hat{\gamma}_{11}{F}_{12}(t)+\hat{\beta}_{1,1}G_{11}(t), \\ \text{if } 0< \alpha_1\leq1,0<\alpha_2\leq2,\qquad
		\\
		\gamma_{10}F_{11}(t)	
		+ \varphi_{1,1}(t-\tau_1)\rho_{1}(t-\tau_1)
	\ast \hat{\gamma}_{11}{F}_{12}(t)+	\hat{\beta}_{1,1}G_{11}(t)
		+\hat{\kappa}_{1,1}G_{12}(t),
		\\
		\text{if } 1< \alpha_1\leq2,0<\alpha_2\leq2,
	\end{array}
	\right.
	\label{rl-delta11eg1}
%\end{aligned}
%\end{eqnarray}
%\begin{eqnarray}
\\%	\begin{aligned}
	\delta_{1,2}(t)=&\left\{
	\begin{array}{lll}
	\hat{\gamma}_{11}	F_{12}(t)\ast \varphi_{1,2}(t-\tau_1)\rho_{1}(t-\tau_1)
		+	\hat{\beta}_{1,2}G_{11}(t),
		 \text{if } 0< \alpha_1\leq1,
		 		 0<\alpha_2\leq2,
		\\
	\hat{\gamma}_{11}	F_{12}(t)\ast \varphi_{1,2}(t-\tau_1)\rho_{1}(t-\tau_1)
	+	\hat{\beta}_{1,2}G_{11}(t) +\hat{\kappa}_{1,2}G_{12}(t), \\
	\text{if } 1< \alpha_1\leq2,0<\alpha_2\leq2,
	\end{array}
	\right.
%	\label{rl-delta12eg1}
	\\
	\delta_{2,i}(t)=&\left\{
	\begin{array}{lll}
	\hat{\gamma}_{22}	F_{21}(t)\ast \varphi_{2,i}(t-\tau_2)\rho_{1}(t-\tau_2)
	+	\hat{\beta}_{2,i}G_{21}(t), \text{if } 0< \alpha_1\leq2,0<\alpha_2\leq1,
		\\
		\hat{\gamma}_{22}	F_{21}(t)\ast \varphi_{2,i}(t-\tau_2)\rho_{1}(t-\tau_2)
+\hat{\beta}_{2,i}G_{21}(t)
		+	\hat{\kappa}_{2,i}G_{22}(t),
	\\	\text{if } 0< \alpha_1\leq2,1<\alpha_2\leq2,
	\end{array}
	\right.	%\label{rl-delta2ieg1}
\end{aligned}
\end{eqnarray}
where
$F_{11},F_{12},F_{21},$ $\varphi_{k,i}(t),\rho_{1}(t)$ are given in \eqref{f1i},
$$
\begin{aligned}
&	G_{11}(t)=	\sum\limits_{m=0}^{\lfloor\frac{t}{\tau_{1}}\rfloor}\hat{\gamma}_{11}^m\eta_1^{\alpha_1(m+1)-1}   E^{m+1}_{\alpha_1,\alpha_1(m+1)}(\gamma_{{11}}\eta_1^{\alpha_1}),
\\	&
	G_{12}(t)=	\sum\limits_{m=0}^{\lfloor\frac{t}{\tau_{1}}\rfloor}\hat{\gamma}_{11}^m\eta_1^{\alpha_1(m+1)-2}   E^{m+1}_{\alpha_1,\alpha_1(m+1)-1}(\gamma_{{11}}\eta_1^{\alpha_1}),
 \\	&	G_{21}(t)=	\sum\limits_{m=0}^{\lfloor\frac{t}{\tau_{2}}\rfloor}\hat{\gamma}_{22}^m\eta_2^{\alpha_2(m+1)-1}  E^{m+1}_{\alpha_2,\alpha_2(m+1)}\big(C\eta_2^{\alpha_2}\big),
\\&
	G_{22}(t)=	\sum\limits_{m=0}^{\lfloor\frac{t}{\tau_{2}}\rfloor}\hat{\gamma}_{22}^m\eta_2^{\alpha_2(m+1)-2}   E^{m+1}_{\alpha_2,\alpha_2(m+1)-1}\big(C\eta_2^{\alpha_2}\big),\eta_k=t-\tau_km,k=1,2,
\end{aligned}
 $$ and $\ast$ is the convolution product of functions.
Hence, we get the analytical solution for the given system \eqref{DRE:eg1} in the sense of the Riemann-Liouville {fractional partial time derivative} as follows:
\begin{eqnarray*}
	\begin{aligned}
		&	u_1(x,t)=\delta_{1,1}(t)+\delta_{1,2}(t)x, \text{and }
		u_2(x,t)=\delta_{2,1}(t)sin(\sqrt{\mu_{2,0}}x)+\delta_{2,2}(t)cos(\sqrt{\mu_{2,0}}x),
	\end{aligned}
\end{eqnarray*}
where $\delta_{k,i}(t),i,k=1,2,$ are given in \eqref{rl-delta11eg1}.
The above-analytical solution of  \eqref{DRE:eg1} in the sense of the  Riemann-Liouville {fractional partial time derivative} obeys the following initial conditions:
\begin{eqnarray*}
	\begin{aligned}
	&	{^{RL}\partial}^{\alpha_1-1}_{t}	[u_1(x,t)]\Big{|}_{t=0}=\hat{\beta}_{1,1}+ \hat{\beta}_{1,2}x, \quad %\text{and }
	{^{RL}\partial}^{\alpha_1-2}_{t}	[u_1(x,t)]\Big{|}_{t=0}=\hat{\kappa}_{1,1}+ \hat{\kappa}_{1,2}x,
\\&	{^{RL}\partial}^{\alpha_2-1}_{t}[u_2(x,t)]\Big{|}_{t=0}= \hat{\beta}_{2,1}sin(\sqrt{\mu_{2,0}}x)
		+\hat{\beta}_{2,2}cos(\sqrt{\mu_{2,0}}x),
		\\& {^{RL}\partial}^{\alpha_2-2}_{t}[u_2(x,t)]\Big{|}_{t=0}= \hat{\kappa}_{2,1}sin(\sqrt{\mu_{2,0}}x)
	+\hat{\kappa}_{2,2}cos(\sqrt{\mu_{2,0}}x),
		\end{aligned}
\end{eqnarray*}
\begin{eqnarray*}
\begin{aligned}
	\text{and the boundary conditions: }		&	u_1(0,t)=\delta_{1,1}(t) %\text{and }
		, u_1(\lambda,t)=\delta_{1,1}(t)+\lambda\delta_{1,2}(t), \\&	
			%\end{aligned}\begin{eqnarray}
		%\begin{aligned}
	u_2(0,t)=\delta_{2,2}(t),	u_2(\lambda,t)=\lambda_1\delta_{2,1}(t) +\lambda_2\delta_{2,2}(t),
	\end{aligned}
\end{eqnarray*}
where $\delta_{k,i}(t),i,k=1,2,$ are given in \eqref{rl-delta11eg1}, $\lambda_1=sin(\lambda\sqrt{\mu_{2,0}})$ and $\lambda_2=cos(\lambda\sqrt{\mu_{2,0}}).$
Note that   the derived solutions  \eqref{form of sol eg1} of the given system   \eqref{DRE:eg1}  along with two  fractional partial time derivatives in the sense of the Caputo and Riemann-Liouville  are valid for   integer values of $\alpha_k=1,$   and  $\alpha_k=2,k=1,2.$
\end{example}

\begin{example}
	\label{eg3}
	Finally, we consider the following  {coupled nonlinear time-fractional DR system} with  linear time delays
	
	\begin{eqnarray}
		\begin{aligned}\label{DRE:eg3}
			{\partial^{\alpha}_{t}} u_k=%\mathcal{R}_1(\mathbf{U},\hat{\mathbf{U}})\equiv
			&\gamma_{k1}u_1
			+ \partial_x[(a_{12}u_2+a_{11}u_1+a_{10})\partial_x u_1
			+(b_{k2}u_2+b_{k1}u_1+b_{k0})\partial_x u_2 ]+\gamma_{k0}
			\\&	+\gamma_{k2}u_2 +
			\hat{\gamma}_{k1}\hat{u}_1+\hat{\gamma}_{k2}\hat{u}_2,\alpha\in(0,2],
\\\hat{u}_k(x,t)=&\phi_k(x,t),t\in[-\tau,0],\tau>0,k=1,2.
		\end{aligned}
	\end{eqnarray}
The discussed Case \ref{poly case} in the previous subsection {confirms} that the  above-system \eqref{DRE:eg3} is invariant under the following vector space
\begin{eqnarray}
	\label{vector space example 3}
	\mathbf{Y}_4= \mathcal{L}\{1,x\} \times \mathcal{L}\{1,x\}.
\end{eqnarray}
Now, we follow the similar procedure  as discussed in Example \ref{eg1}. Thus, the obtained  analytical solution for the given system \eqref{DRE:eg3} with $a_{k2}=-b_{k1},	\gamma_{21}=\hat{\gamma}_{21}=a_{k1}=b_{k2}=0,k=1,2,$ is as follows:
\begin{eqnarray}\label{form of soln eg3}
	\begin{aligned}
		&	u_k(x,t)=\delta_{k,1}(t)+\delta_{k,2}(t)x, k=1,2,
	\end{aligned}
\end{eqnarray}
where $\delta_{k,i}(t),i=1,2$ can be determined using \eqref{vector space example 3} and \eqref{DRE:eg3}. %Here
Thus, we obtain the solution of \eqref{DRE:eg3} with the  {Caputo  fractional partial time derivative} ${\partial^{\alpha_k}_{t}}={^C\partial}_t^{\alpha_k},k=1,2,$ as follows:
\begin{eqnarray}
	\begin{aligned}
				\label{delta11eg3}
		\delta_{1,1}(t)=&\left\{
	\begin{array}{lll}
		\gamma_{10}F_{11}(t)	
		+	\beta_{1,1}F_{12}(t) +F_{13}(t)
			\ast [\big(\hat{\gamma}_{11}\varphi_{1,1}(t-\tau)
		+\hat{\gamma}_{12}\varphi_{2,1}(t-\tau)\big)
		\\\times \rho_{1}(t-\tau)
		+\hat{\gamma}_{12}\hat{\delta}_{2,1}(t-\tau)
		+{\gamma_{12}}\delta_{2,1}(t)],
	\text{if } 0< \alpha\leq1,
		\\
		\gamma_{10}F_{11}(t)	
	+	\beta_{1,1}F_{12}(t) +F_{13}(t)
		\ast [\big(\hat{\gamma}_{11}\varphi_{1,1}(t-\tau)
	+\hat{\gamma}_{12}\varphi_{2,1}(t-\tau)\big) \\\times \rho_{1}(t-\tau)
 	+\hat{\gamma}_{12}\hat{\delta}_{2,1}(t-\tau)
	+{\gamma_{12}}\delta_{2,1}(t)]
	+\kappa_{1,1}F_{14}(t),
	 \text{if } 1< \alpha\leq2,
	\end{array}
	\right.
\\
	\delta_{2,1}(t)=&\left\{
	\begin{array}{lll}
		\gamma_{20}F_{21}(t)	
		+	\beta_{2,1}F_{22}(t) +\hat{\gamma}_{22}F_{23}(t)
			\ast \varphi_{2,1}(t-\tau)\rho_{1}(t-\tau)
		,
		\text{if } 0< \alpha\leq1,
		\\
		\gamma_{20}F_{21}(t)	
		+	\beta_{2,1}F_{22}(t) +\hat{\gamma}_{22}F_{23}(t)
				\ast \varphi_{2,1}(t-\tau)\rho_{1}(t-\tau)
		+	\kappa_{2,1}F_{24}(t),
		\\\text{if } 1< \alpha\leq2,
	\end{array}
	\right.
\\
	\label{delta22eg3}
	\delta_{2,2}(t)=&\left\{
	\begin{array}{lll}
		\hat{\gamma}_{22}F_{23}(t)\ast \varphi_{2,2}(t-\tau)\rho_{1}(t-\tau)
		+\beta_{2,2}F_{22}(t) ,
		\text{if } 0< \alpha\leq1,
		\\
		\hat{\gamma}_{22}F_{23}(t)\ast \varphi_{2,2}(t-\tau)\rho_{1}(t-\tau)
		+\beta_{2,2}F_{22}(t)	+	\kappa_{2,2}F_{24}(t),
		\text{if } 1< \alpha\leq2,
	\end{array}
	\right.\\
	%\label{delta21eg3}
%\end{aligned}\end{eqnarray}
%\begin{eqnarray}
%\begin{aligned}
	\delta_{1,2}(t)=&\left\{
	\begin{array}{lll}
		F_{13}(t)\ast [\rho_{1}(t-\tau)\big(\hat{\gamma}_{11}\varphi_{1,2}(t-\tau)
			+\hat{\gamma}_{12}\varphi_{2,2}(t-\tau)\big)
					+\hat{\gamma}_{12}\hat{\delta}_{2,2}(t-\tau)	\\	+{\gamma_{12}}\delta_{2,2}(t)
									]+	\beta_{1,2}F_{12}(t)  ,
		 \text{if } 0< \alpha\leq1,
		\\
		F_{13}(t)\ast [{\gamma_{12}}\delta_{2,2}(t)+\big(\hat{\gamma}_{11}\varphi_{1,2}(t-\tau)
		+\hat{\gamma}_{12}\varphi_{2,2}(t-\tau)\big)\rho_{1}(t-\tau)
\\	+\hat{\gamma}_{12}\hat{\delta}_{2,2}(t-\tau)
		]
		+		\beta_{1,2}F_{12}(t) +	\kappa_{1,2}F_{14}(t), \text{if } 1< \alpha\leq2,
	\end{array}
	\right.
	\label{delta12eg3}
\end{aligned}
\end{eqnarray}
where $\beta_{k,i}=\delta_{k,i}(0) ,$ $\kappa_{k,i}=
\mathbf{D}_t[\delta_{k,i}(t)]\Big{|}_{t=0},k,i=1,2,$
 \begin{eqnarray}\label{f1ieg3a}
	&\begin{aligned}	F_{k1}(t)= &	\sum\limits_{m=0}^{\lfloor\frac{t}{\tau}\rfloor}\hat{\gamma}_{kk}^m\eta^{\alpha(m+1)}   E^{m+1}_{\alpha,\alpha(m+1)+1}(\gamma_{{kk}}\eta^{\alpha}),
\,	F_{k2}(t)=	\sum\limits_{m=0}^{\lfloor\frac{t}{\tau}\rfloor}\hat{\gamma}_{kk}^m\eta^{\alpha m}    E^{m+1}_{\alpha,\alpha m+1}(\gamma_{{kk}}\eta^{\alpha})
\\
F_{k3}(t)=&\sum\limits_{m=0}^{\lfloor\frac{t}{\tau}\rfloor}\hat{\gamma}_{11}^{m+1}\eta^{\alpha(m+1)-1}    E^{m+1}_{\alpha,\alpha(m+1)}(\gamma_{{11}}\eta^{\alpha}),
\rho_{1}(t)=\left\{\begin{array}{ll}
	0,&\text{if }t\geq0,
	\\
	1&\text{if }t<0,
\end{array} \right. \\F_{k4}(t)=&	\sum\limits_{m=0}^{\lfloor\frac{t}{\tau}\rfloor}\hat{\gamma}_{kk}^m\eta^{\alpha m+1}    E^{m+1}_{\alpha,\alpha m+2}(\gamma_{{kk}}\eta^{\alpha}),
%\end{aligned}
%\end{eqnarray}
%\begin{eqnarray}\label{f1ieg3b}
%\begin{aligned}
 \delta_{k,i}(t)=\varphi_{k,i}(t),t\in[-\tau_{k},0],
  \end{aligned}\end{eqnarray}  $\varphi_{k,i}(t)=\varphi_{k,i}(0),t>0,i,k=1,2, \eta=t-\tau m, $ { and } $\ast$  { is the convolution product of functions.}
Also, we note that solution \eqref{form of soln eg3}  with above-obtained functions  $\delta_{k,i}(t),i=1,2,$  given in  \eqref{delta11eg3}   satisfies
the following initial and  boundary conditions
\begin{eqnarray*}
	\begin{aligned}
	\label{icbc}&	u_k(x,0)= \beta_{k,1}
	+\beta_{k,2}x, 	u_k(0,t)=\delta_{k,1}(t), \\& \mathbf{D}_t [u_k(x,t)]\Big{|}_{t=0}=\kappa_{k,1}
+\kappa_{k,2}x,
\text{ and }
	u_k(\lambda,t)=\delta_{k,1}(t)+\lambda\delta_{k,2}(t),k=1,2.
\end{aligned}
\end{eqnarray*}
Thus, we obtain the functions $\delta_{k_i}(t),i,k=1,2,$ given in \eqref{form of soln eg3} in the sense of {Riemann-Liouville fractional partial time derivative } ${\partial^{\alpha}_{t}}={{^{RL}\partial}}^{\alpha}_{t}$  as follows:
\begin{eqnarray}\begin{aligned}		
	\delta_{1,1}(t)=&\left\{
	\begin{array}{lll}
		\gamma_{10}F_{11}(t)	
		+	\hat{\beta}_{1,1}G_{11}(t) +F_{13}(t)
			\ast [\rho_{1}(t-\tau)\big(\hat{\gamma}_{12}\varphi_{2,1}(t-\tau)
		\\
		+\hat{\gamma}_{11}\varphi_{1,1}(t-\tau)\big)
		+\hat{\gamma}_{12}\hat{\delta}_{2,1}(t-\tau)
	+{\gamma_{12}}\delta_{2,1}(t)]
		,
		\text{if } 0< \alpha\leq1,
		\\
		\gamma_{10}F_{11}(t)	
		+\hat{\beta}_{1,1}G_{11}(t) +	\hat{\kappa}_{1,1}G_{12}(t)+F_{13}(t)
		\ast [\rho_{1}(t-\tau)\big(
		\hat{\gamma}_{12}\varphi_{2,1}(t-\tau)
		\\+ \hat{\gamma}_{11}\varphi_{1,1}(t-\tau)\big)
		+\hat{\gamma}_{12}\hat{\delta}_{2,1}(t-\tau)
			+{\gamma_{12}}\delta_{2,1}(t)]
		, \text{if } 1< \alpha\leq2,
	\end{array}
	\right.	
	\\\delta_{2,1}(t)=&\left\{
	\begin{array}{lll}
		\gamma_{20}F_{21}(t)	
		+\hat{\beta}_{2,1}G_{21}(t) +\hat{\gamma}_{22}F_{23}(t)
		\ast \varphi_{2,1}(t-\tau)\rho_{1}(t-\tau)
		,
		\text{ if } 0< \alpha\leq1,
		\\
		\gamma_{20}F_{21}(t)	
		+	\hat{\beta}_{2,1}G_{21}(t) +\hat{\gamma}_{22}F_{23}(t)
		\ast \varphi_{2,1}(t-\tau)\rho_{1}(t-\tau)
		+	\hat{\kappa}_{2,1}
		G_{22}(t),\\ \text{if } 1< \alpha\leq2,\label{delta21eg3}
	\end{array}
	\right.
		\label{solution-eg3rl}
\\%	\label{delta11eg3rl}
	\delta_{1,2}(t)=&\left\{
	\begin{array}{lll}
		\hat{\beta}_{1,2}G_{11}(t) +F_{13}(t)
		\ast [\big(\hat{\gamma}_{12}\varphi_{2,2}(t-\tau)\
		+\hat{\gamma}_{11}\varphi_{1,2}(t-\tau)\big)\rho_{1}(t-\tau)
	\\	+{\gamma_{12}}\delta_{2,2}(t)
		+\hat{\gamma}_{12}\hat{\delta}_{2,2}(t-\tau)] ,
		\text{ if } 0< \alpha\leq1,
		\\
		\hat{\beta}_{1,2}G_{11}(t) +F_{13}(t)\ast [\big(\hat{\gamma}_{12}\varphi_{2,2}(t-\tau)
	+\hat{\gamma}_{11}\varphi_{1,2}(t-\tau)
		\big)\rho_{1}(t-\tau)
		\\	+{\gamma_{12}}\delta_{2,2}(t)
		+\hat{\gamma}_{12}\hat{\delta}_{2,2}(t-\tau)]+	\hat{\kappa}_{1,2}G_{12}(t),
		\text{if } 1< \alpha\leq2,
	\end{array}
	\right.
	%	\label{delta21eg3rl}
\\
\label{rl solution}
%label{delta22eg3rl}
		\delta_{2,2}(t)=&\left\{
		\begin{array}{lll}
			\hat{\gamma}_{22}F_{23}(t)\ast \varphi_{2,2}(t-\tau)\rho_{1}(t-\tau)
	+
			\hat{\beta}_{2,2}G_{21}(t),  \text{ if } 0< \alpha\leq1,
			\\
			\hat{\gamma}_{22}F_{23}(t)\ast \varphi_{2,2}(t-\tau)\rho_{1}(t-\tau)
			+	\hat{\beta}_{2,2}
			G_{21}(t) +\hat{\kappa}_{2,2}
			G_{22}(t), \text{ if } 1< \alpha\leq2,
		\end{array}
		\right.
	\end{aligned}
\end{eqnarray}
where $\hat{\beta}_{k,i}={^{RL}D_{t}^{\alpha-1}}[\delta_{k,i}(t)]\Big{|}_{t=0},$ $\hat{\kappa}_{k,i}={^{RL}D_{t}^{\alpha-2}}[\delta_{k,i}(t)]\Big{|}_{t=0},k,i=1,2,$  $ F_{kl}(t),l=1,3,$ are given in \eqref{f1ieg3a},
$$
\begin{aligned}
	G_{k1}(t)=&		\sum\limits_{m=1}^{\lfloor\frac{t}{\tau}\rfloor}\hat{\gamma}_{kk}^{m-1}\eta^{\alpha m-1}   E^{m}_{\alpha,\alpha m}(\gamma_{{kk}}\eta^{\alpha})\ \&\
		G_{k2}(t)=	\sum\limits_{m=1}^{\lfloor\frac{t}{\tau}\rfloor}\hat{\gamma}_{kk}^{m-1}\eta^{\alpha m-2}    E^{m}_{\alpha,\alpha m-1}(\gamma_{{kk}}\eta^{\alpha}),
\end{aligned}$$
where $\eta=t-\tau m,k=1,2.$
Also, we note that the solution \eqref{form of soln eg3} of \eqref{DRE:eg3} with the
 above-obtained $\delta_{k_i}(t),i=1,2,$  given in  \eqref{delta21eg3}  obeys the following initial and the  boundary conditions
\begin{eqnarray*}
	\begin{aligned}
	\label{icbcrl}&	{{^{RL}D}^{\alpha-1}_{t}}u_k(x,t)|_{t=0}={\hat{\beta}_{k,1}}
	+{\hat{\beta}_{k,2}}x,
	 &{{^{RL}D}^{\alpha-2}_{t}}u_k(x,t)|_{t=0}={\hat{\kappa}_{k,1}}
	+{\hat{\kappa}_{k,2}}x,
	\\
	&	u_k(0,t)=\delta_{k,1}(t), \text{  and }
	 &	u_k(\lambda,t)=\delta_{k,1}(t)+\lambda\delta_{k,2}(t),k=1,2.
	\end{aligned}
\end{eqnarray*}
Additionally, we observe that the above-obtained solutions  \eqref{form of soln eg3} of the considered system   \eqref{DRE:eg3} with both the Caputo and Riemann-Liouville fractional partial time derivatives are valid for   integer values of $\alpha=1,$   and  $\alpha=2.$
\end{example}
\begin{remark}
 {
It is important to note that exact solutions of the given system \eqref{DREs} under two different fractional partial time derivatives in the sense of the Caputo and Riemann-Liouville  have not been discussed in the literature. Hence the obtained solutions \eqref{form of sol eg1} and \eqref{form of soln eg3} of the given systems \eqref{DRE:eg1} and \eqref{DRE:eg3} in this work are new and significant since they can be expressed  in terms of elementary functions and special functions like the Mittag-Leffler and the Euler-gamma  functions. Also, we observe that in Example \ref{eg1}  and  Example \ref{eg3},  the obtained  solutions in the sense of both the Caputo and Riemann-Liouville fractional partial time derivatives are different for each  $\alpha_k\in(0,2]$ and $\alpha\in(0,2]$, respectively. Additionally, we note that the obtained solutions \eqref{form of sol eg1} and \eqref{form of soln eg3}  are highly dependent on the delay constants $\tau_k>0,k=1,2,$   in Example \ref{eg1} and $\tau>0$  in Example \ref{eg3}, which appear in the corresponding discussed systems \eqref{DRE:eg1} and \eqref{DRE:eg3}, respectively. We believe that the obtained solutions will be useful in the future for studying numerical solutions of discussed coupled nonlinear time-fractional DR systems with time delays through different numerical methods.\\
In addition, we observe that the obtained solutions \eqref{form of sol eg1} and \eqref{form of soln eg3}  of the  given systems \eqref{DRE:eg1} %in Example \ref{eg1}
and \eqref{DRE:eg3} %in Example   \ref{eg3}
	    act in different manners due to different ranges of $ \alpha_k $ in Example \ref{eg1} and $\alpha $  in Example   \ref{eg3}, which are given below.}
	  \begin{itemize}
\item[$\bullet$]
	 {It can be observed that for the coupled nonlinear time-fractional DR system with positive delays $\tau_k$ and $\tau $,  the derived solutions \eqref{form of sol eg1} and \eqref{form of soln eg3} of the considered systems \eqref{DRE:eg1} in Example \ref{eg1}  and \eqref{DRE:eg3}  in Example   \ref{eg3} behave like an exponential-type function due to the involvement of the three-parameter Mittag-Leffler functions when the ranges of $\alpha_k$ and $\alpha$ increase from $0$ to $1$, respectively.}
\item[$\bullet$] 
 {Also, it can be noted that in the sense of both Caputo and Riemann-Liouville fractional-order derivatives,  
the obtained solutions \eqref{form of sol eg1}  for  \eqref{DRE:eg1} with positive delay  $\tau_k$ in Example \ref{eg1} and  \eqref{form of soln eg3} for \eqref{DRE:eg3} with $\tau>0$ in Example \ref{eg3} exhibit a changed form from exponential-type function behaviour to wave-type function behaviour due to the involvement of the three-parameter Mittag-Leffler functions when $\alpha_k,k=1,2,$ and $\alpha$ increase from $1$ to $2$, respectively.
 From this, we can conclude that the system \eqref{DREs} will help to study the intermediate process between diffusion and wave phenomena if $\alpha_k\in(1,2],k=1,2$. Hence the system \eqref{DREs} can be referred to as coupled nonlinear time-fractional diffusion-reaction-wave type equations when $\alpha_k\in(1,2],k=1,2$.	}
\item[$\bullet$]  {	Moreover, in the sense of both Caputo and Riemann-Liouville fractional-order derivatives, the derived solutions \eqref{form of sol eg1}  of  \eqref{DRE:eg1} in Example \ref{eg1}  and  \eqref{form of soln eg3} of \eqref{DRE:eg3} in Example \ref{eg3}  are also valid  for integer-order cases $\alpha_k=\alpha=1$ and $\alpha_k=\alpha=2,$ respectively.}
  	   		 \end{itemize}
\end{remark}
\section{Some extensions of the invariant subspace method}\label{several}
This section presents some extensions of the  invariant subspace method to obtain invariant vector spaces and analytical solutions of the {$2$-component {coupled nonlinear time-fractional systems} of   PDEs} with several time delays involving linear delay terms. Additionally, we present the details of how to find invariant vector spaces and analytical solutions of {multiple-component {coupled nonlinear time-fractional systems} of PDEs} with several time delays using the invariant subspace method.
\subsection{Extension of the invariant subspace method to 2-component coupled system with several time delays}
Let us consider the $2$-component {coupled nonlinear time-fractional system} of   PDEs with several time delays as follows:
\begin{eqnarray}\begin{aligned}\label{1+1:severaldelay}
		\left( 	{\partial^{\alpha_1}_{t}} u_1,
		{\partial^{\alpha_2}_{t}} u_2
		\right) =&{\mathbf{{{P}}}}[\mathbf{U},\hat{\mathbf{U}}]
		,
		\alpha_1,\alpha_2>0,
		\\	u_k(x,t)=&\phi_k(x,t),t\in[-\tau_k,0],k=1,2,
	\end{aligned}
\end{eqnarray}
where $	{\partial^{\alpha_k}_{t}}(\cdot)$ denotes  $\alpha_k$-th order  {fractional partial time derivative}   of  $\mathbf{U}=(u_1,u_2)$ in the sense of either Caputo ${^C\partial^{\alpha_k}_{t}}(\cdot)$ given in \eqref{Caputo} or Riemann-Liouville  ${^{RL}\partial^{\alpha_k}_{t}}(\cdot)$ given in \eqref{RL},
\begin{eqnarray}\label{pk}
	\begin{aligned}
	{\mathbf{{{P}}}}[\mathbf{U},\hat{\mathbf{U}}]
	=&
	\left(
	{\mathcal{{{P}}}_1}[\mathbf{U},\hat{\mathbf{U}}],
	{\mathcal{{{P}}}_2}[\mathbf{U},\hat{\mathbf{U}}]
	\right) ,
	\end{aligned}
\end{eqnarray}
$${\mathcal{{{P}}}_k}[\mathbf{U},\hat{\mathbf{U}}]
=	{\mathcal{{{P}}}_k}\big( x,
u_1,
\partial_x u_1,
\dots,
\partial_x^{m_k} u_1,
u_2,
\partial_x u_2,\dots,	\partial_x^{m_k}  u_2
\big) +\sum\limits_{j=1}^{2}\sum\limits_{i=1}^{p_j} \hat{\eta}^k_{ji}\hat{u}_{ji}$$
is the nonlinear delay differential operator of order $m_k,m_k\in\mathbb{N},$ with $ \partial_x^{j}(\cdot)=\dfrac{\partial^j}{\partial x^j}(\cdot),j=1,2,\dots,m_k,$
$\hat{\mathbf{U}}=(\hat{u}_{11},\hat{u}_{21},\dots,\hat{u}_{1p_1},\hat{u}_{2p_2}),$
$ u_k=u_k(x,t),$
$\hat{u}_{ki}={u}_{k}(x,t-\tau_{ki}),$
$x,\hat{\eta}^k_{ji}\in\mathbb{R},$ $ \tau_{ki},t>0,$  and  $\tau_k=\text{max}\{\tau_{k1},\tau_{k2},\dots,\tau_{kp_k}\},$
$i=1,2,\dots,p_k,k=1,2.$

Note that  the vector nonlinear delay differential operator  ${\mathbf{{{P}}}}[\mathbf{U},\hat{\mathbf{U}}]
=
 (
{\mathcal{{{P}}}_1}[\mathbf{U},\hat{\mathbf{U}}],
{\mathcal{{{P}}}_2}[\mathbf{U},\hat{\mathbf{U}}]
 ) $ admits the Cartesian product space  $\mathbf{Y}_n$ given in \eqref{vector space}   if
either ${\mathbf{{{P}}}}[\mathbf{U},\hat{\mathbf{U}}]\in \mathbf{Y}_n$ for all $\mathbf{U}\in \mathbf{Y}_n,$ or
${\mathcal{{{P}}}_k}[\mathbf{Y}_n]\subseteq Y_{k,n_k},k=1,2.$
%then  $\mathbf{Y}_n$  is called an invariant vector space for ${\mathbf{{{P}}}}[\mathbf{U},\hat{\mathbf{U}}].$
The following theorem provides the generalized separable form of analytical solution for the system \eqref{1+1:severaldelay}  whenever $\mathbf{Y}_n$  is called an invariant vector space  for ${\mathbf{{{P}}}}[\mathbf{U},\hat{\mathbf{U}}].$
\begin{theorem}
	\label{theorem several delays}
	Let  ${\mathbf{{{P}}}}[\mathbf{U},\hat{\mathbf{U}}]
$ given in \eqref{pk} be an invariant under $\mathbf{Y}_n$   given in \eqref{vector space}.
	Then, the  given system \eqref{1+1:severaldelay}
	admits the generalized separable analytical solution
	as follows:
	\begin{eqnarray*}
		\mathbf{U}(x,t)=\left( \sum\limits_{i=1}^{n_1}\delta_{1,i}(t)\psi_{1,i}(x),\sum\limits_{i=1}^{n_2}\delta_{2,i}(t)\psi_{2,i}(x)\right) ,
	\end{eqnarray*}
	such that the functions $\delta_{k,i}(t),i=1,2,\dots,n_k$ satisfy by the following system of fractional-order time derivative delay ODEs
	\begin{eqnarray*}	\begin{aligned}
			D^{\alpha_k}_t{\delta_{k,i}(t)}&=
			\Lambda_{k,i}\Big(
			\delta_{1,1}(t),
			\delta_{2,1}(t),
			\dots, \delta_{1,n_1}(t)
			,
		 	\delta_{2,n_2}(t)
			\Big)
			+	\hat{	\Lambda}_{k,i}\Big(
			\hat{\eta}^1_{11}
			\hat{\delta}_{1,1}(t-\tau_{11}),\\&
			\hat{\eta}^2_{11}
			\hat{\delta}_{1,1}(t-\tau_{11}),
			\hat{\eta}^1_{21}
			\hat{\delta}_{2,1}(t-\tau_{21}),
			 	\hat{\eta}^2_{21}
			\hat{\delta}_{2,1}(t-\tau_{21}),
			\dots, \hat{\eta}^1_{1p_1}
			\hat{\delta}_{1,n_1}(t-\tau_{1p_1}),
			\\&	
			\hat{\eta}^2_{1p_1}
			\hat{\delta}_{1,n_1}(t-\tau_{1p_1}),
			\hat{\eta}^1_{2p_2}
			\hat{\delta}_{2,n_2}(t-\tau_{2p_2}),   \hat{\eta}^2_{2p_2}
			\hat{\delta}_{2,n_2}(t-\tau_{2p_2})\Big),
		\end{aligned}
	\end{eqnarray*}
	for some arbitrary functions $	\Lambda_{k,i}(\cdot)$ and linear functions  $	\hat{\Lambda}_{k,i}(\cdot),$ which can be obtained from the expansion with the basis of $\mathbf{Y}_n.$
\end{theorem}
\begin{proof} Similar to  the proof of Theorem \ref{main result}.
\end{proof}
Next, we present one more extension of the invariant subspace method to obtain solutions for the multiple-component {coupled nonlinear time-fractional systems} of   PDEs with several time delays in the following subsection.
\subsection{One more extension of the invariant subspace method to multiple-component coupled system with several time delays}
In this subsection, we consider the multiple-component {coupled nonlinear time-fractional system} of   PDEs with several time delays involving linear delay terms as follows:
\begin{eqnarray}\begin{aligned}\label{multi:severaldelay}
		&\left( 	{\partial^{\alpha_1}_{t}} u_1,
		\dots,	{\partial^{\alpha_q}_{t}} u_q
		\right)
		={\Phi}[\mathbf{U},\hat{\mathbf{U}}]
		,\alpha_k>0,
		\\&	u_k(x,t)=\phi_k(x,t),t\in[-\tau_k,0],k=1,2,\dots,q,
	\end{aligned}
\end{eqnarray}
where $	{\partial^{\alpha_k}_{t}}(\cdot),1,2,\dots,q,  $ represents  $\alpha_k$-th order  { fractional partial time derivative}   of  $\mathbf{U}=(u_1,u_2,\dots,u_q), $ in the sense of either Caputo ${^C\partial^{\alpha_k}_{t}}(\cdot)$ given in \eqref{Caputo} or  Riemann-Liouville  ${^{RL}\partial^{\alpha_k}_{t}}(\cdot)$ given in \eqref{RL},  ${\Phi}[\mathbf{U},\hat{\mathbf{U}}]
=
 (
\varPhi_1[\mathbf{U},\hat{\mathbf{U}}],
\dots,
\varPhi_q[\mathbf{U},\hat{\mathbf{U}}]
 ),$
\begin{eqnarray*}\label{varphik}
	\begin{aligned}
		\varPhi_k[\mathbf{U},\hat{\mathbf{U}}]
		=&	\varPhi_k\big( x,
		u_1,
		\partial_x u_1,
		\dots,
		\partial_x^{m_k} u_1,\dots,u_q,
	 \partial_x u_q,
		\dots,\partial_x^{m_k} u_q
		\big)  +\sum\limits_{l=1}^{q}\sum\limits_{i=1}^{p_l} \hat{\eta}^k_{li}\hat{u}_{li},
	\end{aligned}
\end{eqnarray*}
is the nonlinear delay differential operator of order $m_k,m_k\in\mathbb{N},$ with $\hat{\mathbf{U}}=(\hat{u}_{11},\hat{u}_{21},\dots,\hat{u}_{q1},\dots,\hat{u}_{1p_1},\hat{u}_{2p_2},\dots,\hat{u}_{qp_q}),$
$\hat{u}_{ki}={u}_i(x,t-\tau_{ki}), $ $u_k=u_k(x,t),$ $\hat{\eta}^k_{ji},x\in\mathbb{R}, $ $\tau_{ki},t>0,$ $i=1,2,\dots,p_k,$ $ \partial_x^{j}(\cdot)=\dfrac{\partial^j}{\partial x^j}(\cdot),j=1,2,\dots,m_k,$
and  $\tau_k=\text{max}\{\tau_{k1},\tau_{k2},\dots,\tau_{kp_k}\},$
$k=1,2,\dots,q.$

Now, we consider a Cartesian product space of $q$-component vector spaces \begin{eqnarray}\label{multi- vector space}
	\hat{\mathbf{Y}}_n=\hat{Y}_{1,n_1}\times\hat{Y}_{2,n_2}\times\cdots\times\hat{Y}_{q,n_q},
\end{eqnarray}
where $\hat{Y}_{k,n_k}$ is the $n_k$-dimensional component vector space, that is
$
	Y_{k,n_k}=\mathcal{L}
	\{\psi_{k,1}(x), \dots,\psi_{k,n_k}(x)\},k=1,\dots,q.
$
Note that $\{\psi_{k,1}(x),\psi_{k,2}(x),\dots,\psi_{k,n_k}(x)\}$ is the set of $n_k$-linearly independent solutions of the system of  linear homogeneous integer-order ODEs
\begin{eqnarray*}
	\begin{aligned}
		\hat{H}_k[\lambda(x)]	:=&{\mathbf{D}^{n_k}_x\lambda(x)}+\mu_{k,n_k-1}{\mathbf{D}^{n_k-1}_x\lambda(x)} +\dots+\mu_{k,0}\lambda(x)=0,
	\end{aligned}
\end{eqnarray*}
where $\mu_{k,j}\in\mathbb{R},j-1,2,\dots,n_k-1, k=1,2,\dots,q.$
It is known that the $n$-dimensional Cartesian product space $\hat{\mathbf{Y}}_n$ given in \eqref{multi- vector space} is an invariant vector space under the vector nonlinear delay differential operator ${\Phi}[\mathbf{U},\hat{\mathbf{U}}]
=
\left(
\varPhi_1[\mathbf{U},\hat{\mathbf{U}}],
\dots,
\varPhi_q[\mathbf{U},\hat{\mathbf{U}}]
\right)$ if either one of the following  conditions holds:
$${\Phi}[\hat{\mathbf{Y}}_n,\hat{\mathbf{Y}}_n]\subseteq \hat{\mathbf{Y}}_n,\text{ or  }
\varPhi_k[\mathbf{U},\hat{\mathbf{U}}]\in \hat{Y}_{k,n_k},k=1,2,\dots,q, \, \forall \mathbf{U}\in \hat{\mathbf{Y}}_n.$$
The following theorem shows that we can obtain the generalized separable analytical solutions  for a given {coupled nonlinear time-fractional system} \eqref{multi:severaldelay} using the invariant subspace method.
\begin{theorem}
	\label{{multi:severaldelay}:thm}
	Suppose that  ${\Phi}[\mathbf{U},\hat{\mathbf{U}}]
	=
	\left(
	\varPhi_1[\mathbf{U},\hat{\mathbf{U}}],
	\dots,
	\varPhi_q[\mathbf{U},\hat{\mathbf{U}}]
	\right)$  given in \eqref{multi:severaldelay} admits the invariant vector space $\hat{\mathbf{Y}}_n$   given in \eqref{multi- vector space}.
	Then, the multiple-component time-fractional system	\eqref{multi:severaldelay} with several time delays
	has a  generalized separable analytical solution
	of the form
	\begin{eqnarray*}
		\begin{aligned}
			&	\mathbf{U}(x,t)=(u_1,u_2,\dots,u_q),
			 u_k:=u_k(x,t)
			= \sum\limits_{i=1}^{n_k}\delta_{k,i}(t)\psi_{k,i}(x),k=1,2,\dots,q,
		\end{aligned}
	\end{eqnarray*}
	such that the functions $\delta_{k,i}(t),i=1,2,\dots,n_k$ are determined by the following system of fractional-order time derivative delay ODEs
	\begin{eqnarray*}	\begin{aligned}
			D^{\alpha_k}_t{\delta_{k,i}(t)}&=
			\Lambda_{k,i}\Big(
			\delta_{1,1}(t),
			\dots,	\delta_{q,1}(t),
			\dots, \delta_{1,n_1}(t)
			,
			 \dots,\delta_{q,n_q}(t)
			\Big)
			+	\hat{	\Lambda}_{k,i}\Big(
			\hat{\eta}^1_{11}
			\hat{\delta}_{1,1}(t-\tau_{11}),\dots
			\\&		\hat{\eta}^q_{11}
			\hat{\delta}_{1,1}(t-\tau_{11}),
			\dots	\hat{\eta}^1_{q1}
			\hat{\delta}_{2,1}(t-\tau_{q1}),
			\dots 		\hat{\eta}^q_{q1}
			\hat{\delta}_{q,1}(t-\tau_{q1}),
			\dots,  \hat{\eta}^q_{1p_1}
			\hat{\delta}_{1,n_1}(t-\tau_{1p_1}),\\&
			\dots,	\hat{\eta}^q_{qp_q}
			\hat{\delta}_{q,n_q}(t-\tau_{qp_q})\Big),
		\end{aligned}
	\end{eqnarray*}
	for some arbitrary functions $	\Lambda_{k,i}(\cdot)$ and linear functions  $	\hat{\Lambda}_{k,i}(\cdot),i=1,2,\dots,n_k,$ $k=1,2,\dots,q,$ which can be obtained from the expansion with the basis of $\mathbf{Y}_n.$
\end{theorem}
\begin{proof} Similar to  the proof of Theorem \ref{main result}.
\end{proof}
\section{Concluding remarks}\label{sec4}
In this work,  we systematically discussed the application  of the invariant subspace method for {obtaining} analytical solutions of  given {coupled  nonlinear  time-fractional system of   PDEs} \eqref{1+1:delay} with  time delays.
Also, we  presented  a systematic way to obtain different dimensional invariant vector spaces $\mathbf{Y}_n$ for the given {system} \eqref{1+1:delay} using the system of linear homogeneous integer-order ODEs given in \eqref{lhodes}. Additionally, we observed that using the obtained invariant vector spaces $\mathbf{Y}_n$, we can reduce the {given system} \eqref{1+1:delay}    into  the system of fractional-order time derivative delay ODEs.
 More specifically, we explicitly illustrated the effectiveness of
 this method by finding analytical solutions for the {IBVPs} of the
{coupled nonlinear time-fractional DR system} \eqref{DREs} with
 linear time delays
 under
 the two {distinct} fractional derivatives that are (a) the Riemann-Liouville fractional partial time
 derivative and (b) the Caputo fractional partial time derivative with  different orders  $\alpha_1$ and $\alpha_2,$  where $\alpha_1,\alpha_2\in(0,2].$
 %  We observe that invariant subspace method can be  efficiently applied to obtain analytical solutions for coupled system of DREs  \eqref{DREs} with multiple delays. Analytical solutions obtained  were different  under the two definitions say Caputo and Riemann-Liouville time-fractional partial derivatives for $\alpha_k\in(0,2],k=1,2.$   The major difference in the obtained analytical solutions for a given coupled system of DREs  \eqref{DREs} with multiple delays is due to the difference in the initial conditions with respect the Caputo and Riemann-Liouville time-fractional partial derivatives for $\alpha_k\in(0,2],k=1,2.$  Caputo time-fractional partial derivative in \eqref{DREs} provides solutions with initial conditions in integer-order derivatives of $u_k(x,t),k=1,2,$ where as  Riemann-Liouville time-fractional partial derivative in \eqref{DREs} gives solutions having initial conditions with time-fractional derivative of $u_k(x,t),k=1,2.$ This study shows the potential of invariant subspace method to derive analytical solutions for various nonlinear fractional coupled systems of PDEs with and without time delays under various fractional derivatives in the sense of Prabhakar derivative, Hilfer derivative, regularized Prabhakar derivative, and Weyl derivative.
    {Additionally, we have shown some extensions of the   invariant subspace method to obtain analytical solutions of $2$-component  and multiple-component {coupled nonlinear time-fractional systems of PDEs} with several time delays as given in \eqref{1+1:severaldelay} and \eqref{multi:severaldelay}, respectively.}

  From this work, we observed the following points for applying the invariant subspace method to {the coupled  nonlinear  time-fractional system of PDEs} \eqref{1+1:delay} with time delays.
  \begin{itemize}
  	\item[$\bullet$]Solving the {coupled  nonlinear  time-fractional system of PDEs} \eqref{1+1:delay} with time delays is very complicated due to the involvement of delay terms and  singular kernel fractional partial time derivatives. However, the invariant subspace method is successfully applied to obtain the analytical solutions for {the given system} \eqref{1+1:delay}  under two fractional-order derivatives that are (a) the Riemann-Liouville   and (b) the Caputo  fractional partial time derivatives.
  	\item[$\bullet$]The derived analytical solutions for the {coupled nonlinear time-fractional DR system} \eqref{DREs} with multiple time delays in the sense of Caputo and Riemann-Liouville are different due to initial conditions involving integer-order derivatives in the sense of Caputo and fractional-order initial conditions involved in the sense of Riemann-Liouville.
  	\item[$\bullet$] Also, this method is applicable to higher-dimensional scalar and {coupled  nonlinear  fractional systems} of  PDEs with and without time delay terms under different kinds of fractional-order derivatives.
  \end{itemize}
  Through this study, we observed that the invariant subspace method is an essential mathematical technique for finding analytical solutions of many significant  {coupled  nonlinear  time-fractional systems} of PDEs with time delays under different  types of  fractional  derivatives. Also, we firmly believe that in future, this method will be very beneficial in finding generalized separable analytical solutions for complex {coupled nonlinear fractional  systems} of PDEs with time delays, which arise in potential applications in diverse areas of science and engineering.

%\section*{Acknowledgement} The first author (K.S.P.) would like to thank the International Mathematical Union (IMU), Germany,  for providing financial support in the form of IMU Breakout Graduate fellowship-2023 (IMU-BGF-2023-06). Another author (M.L.) is supported by a Department of Science and Technology, India-SERB National Science Chair position (NSC/2020/000029).
%\newpage
%	\section*{Declarations}
	
%	Some journals require declarations to be submitted in a standardised format. Please check the Instructions for Authors of the journal to which you are submitting to see if you need to complete this section. If yes, your manuscript must contain the following sections under the heading `Declarations':
	
%	\begin{itemize}
%		\item Funding
%		\item Conflict of interest/Competing interests (check journal-specific guidelines for which heading to use)
%		\item Ethics approval
%		\item Consent to participate
%		\item Consent for publication
%		\item Availability of data and materials
%		\item Code availability
%		\item Authors' contributions
%	\end{itemize}
\section*{Acknowledgments}
The first author (K.S.P.) is thankful for the financial support in the form of International Mathematical Union Breakout Graduate fellowship-2023 (IMU-BGF-2023-06) provided by the  IMU, Germany. Another author (M.L.) is supported by a Department of Science and Technology, India-SERB National Science Chair position (NSC/2020/000029).
	\section*{Declarations}
	\textbf{Conflict of interest: }
The authors declare that they have no conflict of interest.
\\
\textbf{Availability of data and materials: } Not applicable.
	
%	\noindent
%	If any of the sections are not relevant to your manuscript, please include the heading and write `Not applicable' for that section.
	
%\textbf{Ethical Approval:} Not applicable.	\\
%\textbf{Availability of data and materials:} Not applicable.\\
%\textbf{Author contributions:}All authors wrote the main manuscript text, and K.S.P. and P.P. prepared figures. P.P. and M.L. supervised the work.All the authors validated the formal analysis and reviewed the manuscript.

\end{document}